\documentclass[12pt]{amsart}
\usepackage[english]{babel}
\usepackage[all]{xy}
 \usepackage[latin1]{inputenc}
  \usepackage[T1]{fontenc}
  \usepackage{amsmath,amssymb}
  \usepackage{amsmath}
  \usepackage{amsthm}

\usepackage{amssymb}
\usepackage{mathrsfs}
\usepackage{enumerate}
\usepackage{tikz}
\usetikzlibrary{shapes.geometric}
\usepackage{physics}
\usepackage{dsfont}
\newtheorem{theorem}{Theorem}[section]
\newtheorem{problem}[theorem]{Problem}
\newtheorem{lemma}[theorem]{Lemma}
\newtheorem{proposition}[theorem]{Proposition}

\newtheorem{corollary}[theorem]{Corollary}
\newtheorem{definition}[theorem]{Definition}

\theoremstyle{definition}
\newtheorem{remark}[theorem]{Remark}

\usepackage{tikz}
\usetikzlibrary{arrows,shapes,positioning,shadows,trees,matrix}
\usepackage{xcolor}

\title[]{Schur analysis over the unit spectral ball}

\author[D. Alpay]{Daniel Alpay}
\address{(DA) Schmid College of Science and Technology \\
Chapman University\\
One University Drive
Orange, California 92866\\
USA}
\email{alpay@chapman.edu}

\author[I. Cho]{Ilwoo Cho}
\address{(IC) Department of Mathematics and Statistics \\
Saint Ambrose University \\
508 W. Locust St.
Davenport, IA 52803\\
USA}
\email{choilwoo@sau.edu}

\begin{document}

\begin{abstract}
  We begin a study of Schur analysis when the variable is now a matrix rather than a complex number.
  We define the corresponding Hardy space, Schur multipliers and their realizations, and interpolation. Possible applications of the present work
  include matrices of quaternions, matrices of split quaternions,  and other algebras of hypercomplex numbers.
\end{abstract}

\maketitle

\tableofcontents


\section{Introduction}
\setcounter{equation}{0}
The classical theory of Hilbert spaces of  functions analytic in a neighborhood of the origin, i.e., with elements of the form
\begin{equation}
  \label{newform}
  f(z)=\sum_{n=0}^\infty z^n f_n
\end{equation}
where the $f_n\in\mathbb C$ satisfy
\begin{equation}
  \label{gamma-n}
\sum_{n=0}^\infty \gamma_n |f_n|^2<\infty
  \end{equation}
  encompasses spaces such as the Hardy space, the Fock space, the Bergman space, and many others. In \eqref{gamma-n} the numbers $\gamma_n$
  (the weights) are strictly positive numbers such that
  \[
    R=\sqrt{\liminf_{n\rightarrow\infty}\gamma_n^{1/n}}>0.
  \]
  From
  \[
    |\sum_{n=0}^\infty z^nf_n|^2\le\left(\sum_{n=0}^\infty \frac{|z|^{2n}}{\gamma_n}\right)\left(\sum_{n=0}^\infty \gamma_n |f_n|^2\right)
    \]
    the functions are analytic for $z$ such that
\[
  |z|^2< \frac{1}{\limsup_{n\rightarrow\infty}\dfrac{1}{\gamma_n^{1/n}}}=\liminf_{n\rightarrow\infty}\gamma_n^{1/n}=R^2
\]
    ($R=\infty$ is allowed; one then has a space of entire functions). Furthermore the corresponding reproducing kernel is
\begin{equation}
\label{rk-gamma}  
K(z,w)=\sum_{n=0}^\infty\frac{z^n\overline{w}^n}{\gamma_n}.
\end{equation}
One can also allow some of the weights $\gamma_n$ to vanish, as for instance for the Dirichlet kernel
\[
  -\ln(1-z\overline{w})=\sum_{n=1}^\infty \frac{z^n\overline{w}^n}{n},
\]
where $z,w$ belong to the open unit disk $\mathbb D$. Then there is no constant term in the expansions \eqref{gamma-n} and \eqref{rk-gamma}.\smallskip

The theory of such spaces can be extended to the case where the coefficients $f_n$, and possibly the weights $\gamma_n$, are matrices; see e.g., \cite{ab6}. In the present paper we consider the case where, in \eqref{newform}, not only $f_n$ but also $z$ is replaced by a matrix.
Thus we consider expressions of the form
\begin{equation}
 F= F(Z)=\sum_{n=0}^\infty Z^nF_n
  \label{FZequal}
  \end{equation}
where\footnote{We will freely use both notations $F$ and $F(Z)$.} $Z\in\mathbb C^{p\times p}$ and the $F_n\in\mathbb C^{p\times p}$, or more generally belong to a Hilbert space which is also a right and left
$\mathbb C^{p\times p}$-module. When $p=1$ it may be that \eqref{FZequal} converges only for $Z=z=0$. When $p>1$, there are always non-zero nilpotent matrices and the set of convergence of \eqref{FZequal} is not reduced to a point,
but it does not contain necessarily an open set (i.e. may have an empty interior), as follows from Proposition \ref{prop2-3} below. See Corollary \ref{empty-interior}.\smallskip

Expressions of the form \eqref{FZequal} are not closed under pointwise product, and following \cite{MR51:583} we
define on monomials $Z^nA$ and $Z^mB$ the convolution (or Cauchy; see \cite{MR51:583}) product as follows:
Let $n,m\in\mathbb N_0$ and $A,B\in\mathbb C^{p\times p}$. One sets

    \begin{equation}
      (Z^nA)\star(Z^mB)=Z^{n+m}AB,
      \label{cauchy-prod}
    \end{equation}
    extended by linearity to the linear span of the monomials. This is called the Cauchy product and appears 
    in various places in non-commutative settings; see \cite{MR51:583} and, for the quaternionic setting, see \cite{MR2752913} for slice functions and
    \cite{MR618518} for hypercomplex regular functions. \smallskip

    Among other related non-commutating settings let us mention for example the calculus on diagonals developed in \cite{MR92g:94002,MR93b:47027}
    (see \cite{zbMATH03009786} for a predecessor to this calculus).\smallskip

  By specializing to subspaces of $\mathbb C^{p\times p}$ we obtain known cases such as the quaternions and the split quaternions, for $p=2$ and $Z,W$ of the form
  \begin{equation}
    \label{quater}
\begin{pmatrix}z_1&-z_2\\ \overline{z_2}&\overline{z_1}\end{pmatrix} \quad{\rm and}\quad \begin{pmatrix}z_1&z_2\\ \overline{z_2}&\overline{z_1}\end{pmatrix} 
\end{equation}
respectively, and matrices of these for $p$ even, say $p=2h$ ($h\in\mathbb N$) and $z_1$ and $z_2$ being elements of $\mathbb C^{h\times h}$ in \eqref{quater}.
Bicomplex numbers will correspond to matrices of the form
\[
\begin{pmatrix}z_1&-z_2\\ {z_2}&{z_1}\end{pmatrix},
\]
and hyperbolic numbers correspond to the case where $z_1$ and $z_2$ are real in this latter expression. Here too, our approach allows to study matrices
of bicomplex and hyperbolic numbers.\smallskip

The paper consists of seven sections, of which this introduction is the first. In the second section we give first properties of the power series of the form \eqref{FZequal}, and examples such as the counterparts of the Fock space, the Hardy space, the Wiener algebra and rational functions in the present setting.
In the third section we study in greater details the Hardy space. We discuss a version of the
Beurling-Lax theorem in the present setting and discuss the notion of Blaschke factor here. In Section \ref{sec-4} we solve a Nevanlinna-Pick interpolation problem
in the Hardy space. Schur multipliers and their coisometric
realizations are studied in Section \ref{sec-5} using two approaches. The first uses the theory of linear isometry relations and the second is based on operator ranges. In
particular we prove a version of Leech's factorization theorem, needed in the arguments. The notion of Carath\'eodory multipliers is studied
in Section \ref{Cara} with different methods; we reduce the problem to the case of a complex variable. In the last section we discuss various future directions of research and links with hypercomplex analysis.\smallskip

Throughout the paper $B(0,R)$ denotes the open unit disk in the complex plane centered a the origin, and with radius $R$. The spectral
radius of the matrix $A\in\mathbb C^{p\times p}$ will be denoted by $\rho(A)$.

\section{Generalities and first examples}
\setcounter{equation}{0}
\subsection{Preliminary results}

We first introduce:

\begin{definition}
  \label{defhr}
  Given $R>0$,   the ring $\mathcal H_R$ consists of the power series of the form
\begin{equation}
F(Z)=\sum_{n=0}^\infty Z^nF_n, \quad where\,\, Z,F_0,F_1,\ldots\in\mathbb C^{p\times p}
\label{aug23}
\end{equation}
and the latter are such that
\begin{equation}
\label{ineq-23}
\limsup_{n\to\infty} \|F_n\|^{\frac{1}{n}}\le \frac{1}{R}.
\end{equation}
\end{definition}

Note that \eqref{aug23} converges for $Z$ such that $\rho(Z)<R$.

\begin{proposition}
  \label{empty}
For any $A\in\mathbb C^{p\times p}$ with $\rho(A)<R$, the series 
\begin{equation}
F(A)=\sum_{n=0}^\infty A^nF_n
\label{aug23aa}
\end{equation}
converges absolutely.
\end{proposition}
\begin{proof}
We have 
\begin{equation}
  \label{rtyu}
\limsup_{n\to\infty}\|A^n F_n\|^{\frac{1}{n}}\le\lim_{n\to\infty}\|A^n\|^{\frac{1}{n}}\cdot \limsup_{n\to\infty}\|F_n\|^{\frac{1}{n}}
\le \frac{\rho(A)}{R}<1,
\end{equation}
and so the series with running term $\|A^nF_n\|$ converges.
\end{proof}

\begin{proposition}
  \label{prop2-3}
In the above notation, assume
\begin{equation}
  \label{abc1212}
    \limsup_{n\to\infty} \|F_n\|^{\frac{1}{n}}=\infty.
  \end{equation}
  Then,  the set of convergence of \eqref{aug23} does not contain invertible matrices.
  \end{proposition}

  \begin{proof}
Let $A\in\mathbb C^{p\times p}$ be invertible. Then,
  \[
F_n=A^{-n}A^nF_n
\]
and so
\[
\|F_n\|^{1/n}\le \|A^{-n}\|^{1/n}
\|A^nF_n\|^{1/n}
\]
Since $A$ is invertible, it holds that $\lim \|A^{-n}\|^{1/n}=\rho(A^{-1})>0$ and we have
\[
\limsup_{n\to\infty} \|F_n\|^{1/n}=\infty\quad\Longrightarrow\quad \limsup_{n\to\infty} \|A^nF_n\|^{1/n}=\infty
\]
and \eqref{aug23aa} cannot converge.
\end{proof}

  \begin{corollary}
\label{empty-interior}
    If \eqref{abc1212} holds, the set where the power series converges has an empty interior.
    \end{corollary}
  \begin{proof}
This follows from the fact that the set of invertible matrices is dense and open in $\mathbb C^{p\times p}$.
    \end{proof}
\begin{remark}
When $p>1$ besides the case where $A$ is nilpotent, one has examples where $AF_n=0$ for all $n\ge 1$; then $F(A)=F_0$ exists even if \eqref{abc1212} is in force.
  \end{remark}

\begin{proposition}
For any $F,G\in\mathcal H_R$ and $A$ with $\rho(A)<R$,
\begin{equation}
  (F+G)(A)=F(A)+G(A)\quad\mbox{and}\quad (F\star G)(A)=\sum_{n=0}^\infty A^nF(A)G_n,
\label{aug23a}
\end{equation}
where $G=\sum_{n=0}^\infty Z^nG_n$ with $G_0,G_1,\ldots\in\mathbb C^{p\times p}$ and
where $\star$ denotes the Cauchy product \eqref{cauchy-prod}. In particular, $(F\star G)(A)=0$ whenever $F(A)=0$. 
\label{R:eval}
\end{proposition}

\begin{proof}
  We only prove the second claim.
  We have
  \[
  \begin{split}
    F\star G&=F\star\left(\sum_{n=0}^\infty Z^nG_n\right)\\
    &=\sum_{n=0}^\infty (F\star Z^n)G_n\\
    &=\sum_{n=0}^\infty Z^nFG_n,
    \end{split}
  \]
  where the exchange of sum and $\star$-product is justified since
  \[
    \begin{split}
      \|\sum_{n=0}^\infty Z^nFG_n\|\le\sum_{n=0}^\infty \|F\|\cdot\|Z^n\|\cdot\|G_n\|,
    \end{split}
    \]
    and using $\rho(Z)<R$ and \eqref{ineq-23} for $G$. This ends the proof since
  \[
(Z^nFG_n)(A)=A^nF(A)G_n.
  \]
  \end{proof}

  \begin{definition}
    Let $R>0$ and let $F(Z)=\sum_{n=0}^\infty Z^nF_n\in\mathcal H_R$. We define
\begin{equation}
\label{aug23bbbb}
F(zI_p)=\sum_{n=0}^\infty z^nF_n,\quad |z|<R.
\end{equation}
    \end{definition}

    \begin{lemma}
      Let $F,G\in\mathcal H_R$. It holds that
  \begin{equation}
(F\star G)(zI_p)=F(zI_p)G(zI_p), \quad |z|<R.
    \end{equation}
\end{lemma}

\begin{proof}
  We have
\[
  \begin{split}
    \left(\left(\sum_{n=0}^\infty Z^nF_n\right)\star\left(\sum_{n=0}^\infty Z^nG_n\right)\right)(zI_p)&=\left(\sum_{n=0}^\infty Z^n\left(\sum_{k=0}^nF_kG_{n-k}\right)\right)(zI_p)\\
    &=    \sum_{n=0}^\infty z^n\left(\sum_{k=0}^nF_kG_{n-k}\right)\\
    &=\left(\sum_{n=0}^\infty z^nF_n\right)\left(\sum_{n=0}^\infty z^nG_n\right)\\
      &=F(zI_p)G(zI_p),
    \end{split}
  \]
  where the various exchanges of sums hold since $|z|<R$.
\end{proof}

\begin{theorem}
  Let $R>0$ and let $F=\sum_{n=0}^\infty Z^nF_n\in\mathcal H_R$. Let
  \[
F(zI_p)=\left(f_{ij}(z)\right)_{i,j=1}^p
    \]
and let $A\in\mathbb C^{p\times p}$ with spectrum inside $B(0,R)$. Then,
  \begin{equation}
    F(A)=\frac{1}{2\pi i}\int_\gamma (zI_p-A)^{-1}F(zI_p)dz
  \end{equation}
  where $\gamma$ is a closed Jordan curve inside $B(0,R)$, and which encloses the spectrum of $A$.
\end{theorem}

\begin{proof}
  With $F_n=(f_{ij}^{(n)})_{i,j=1}^p$, and with
  \[
f_{ij}(z)=\sum_{n=0}^\infty z^nf_{ij}^{(n)},
    \]
we have
  \[
    F(zI_p)=\sum_{n=0}^\infty z^nF_n=\left(f_{ij}(z)\right)_{i,j=1}^p.
    \]

With $E_{ij}$, $i,j=1,\ldots, p$ denoting the standard basis in $\mathbb C^{p\times p}$ we further have
  \[
    \begin{split}
      \sum_{n=0}^\infty A^nF_n&=\sum_{n=0}^\infty A^n\sum_{i,j=1}^pf_{ij}^{(n)}E_{ij}\\
      &=\sum_{i,j=1}^p\left(\underbrace{\sum_{n=0}^\infty A^nf_{ij}^{(n)}}_{f_{ij}(A)}\right)E_{ij}\\
      &=\sum_{i,j=1}^p \left(\frac{1}{2\pi i}\int_{\gamma}(zI_p-A)^{-1}f_{ij}(z)dz\right)E_{ij}dz\\
      &=        \frac{1}{2\pi i}\int_{\gamma}(zI_p-A)^{-1}\left(\sum_{i,j=1}^pf_{ij}(z)E_{ij}\right)dz\\
            &=        \frac{1}{2\pi i}\int_{\gamma}(zI_p-A)^{-1}F(zI_p)dz.
    \end{split}
    \]
  \end{proof}
\begin{theorem}
  \label{thmeasure}
  Assume that there is a summable positive measure $\mu$ of the form 
  \[
d\mu(z)=m(|z|)dxdy,\quad |r|<R
\]
such that
\begin{equation}
  \label{moments}
 2\pi\int_0^R r^{2n+1}m(r)dr=\gamma_n,\quad n=0,1,\ldots
\end{equation}
Then for $F=\sum_{n=0}^\infty Z^nF_n$ such that \eqref{ineq-23} is in force, it holds that
\begin{equation}
  \int_{|z|<R}(F(zI_p))^*F(zI_p)d\mu(z)=\sum_{n=0}^\infty \gamma_nF_n^*F_n.
  \label{inner-prod-matrix}
\end{equation}
  \end{theorem}

  \begin{proof}
The power series $F(Z)$, and in particular $F(zI_p)$, converges in $|z|<R$ in view of \eqref{ineq-23}. To check \eqref{inner-prod-matrix}, let
    \[
f_N(z)=\sum_{n,m=0}^NF_n^*F_m\overline{z}^nz^m,
\]
By Cauchy inequality, we have that, for $|z|\le r<R$,

\[
  \begin{split}
    |f_N(z)|&\le\sum_{n,m=0}^N |z|^n|z|^m\|F_n\|\|F_m\|\\
    &=\left(\sum_{n=0}^N\frac{|z|^n}{\sqrt{\gamma_n}}\sqrt{\gamma_n}\|F_n\|\right)^2\\
    &\le\left(\sum_{n=0}^N\frac{|z|^{2n}}{\gamma_n}\right)\left(\sum_{n=0}^N\gamma_n\|F_n\|^2\right)\\
    &\le\left(\sum_{n=0}^\infty\frac{r^{2n}}{\gamma_n}\right)\left(\sum_{n=0}^\infty\gamma_n\|F_n\|^2\right).
  \end{split}
\]
Since $r<R$ we have that $\sum_{n=0}^\infty\frac{r^{2n}}{\gamma_n}<\infty$, and we can apply the dominated convergence theorem since $d\mu$ is
assumed to be summable. We then obtain
    \[
      \begin{split}
        \int_{|z|<r}(F(zI_p))^*F(zI_p)d\mu(z)&=\sum_{n,m=0}^\infty F_n^*F_m\iint_{|z|<r}\overline{z}^nz^md\mu(z)\\
        &=\sum_{n,m=0}^\infty F_n^*F_m\int_{\rho=0}^r\int_{\theta=0}^{2\pi} \rho^{n+m}e^{i\theta(n-m)}m(\rho)\rho d\rho d\theta\\
        &=\sum_{n=0}^\infty F_n^*F_n\left(2\pi\int_0^r\rho^{2n+1}m(\rho)d\rho\right)
%
              \end{split}
            \]
 where the integrals are finite thanks to \eqref{moments}.  The monotone convergence theorem allows to let $r$ tend to $R$ and obtain the result using here too \eqref{moments}.
    \end{proof}

\begin{remark}
The previous result can be adapted to other cases, such as Dirichlet, where derivatives appear, or to the case of the Hardy space.
Another case where the previous result can be adapted is the time-varying setting as developed in \cite{MR93b:47027, dede}. Let
\[
  U=\sum_{n=0}^\infty \mathcal Z^nD_n
\]
where $\mathcal Z$ is now the unitary shift  from $\ell_2(\mathbb Z,\mathbb C)$ into itself and $D_0,D_1,\ldots$ are diagonal operators. The Zadeh transform $U(z)$ of $U$ is
\begin{equation}
  \label{zadeh}
  \sum_{n=0}^\infty z^n\mathcal Z^nD_n,\quad z\in\mathbb C
  \end{equation}
  (see \cite{MR1878952,MR2003f:47021}, and the unpublished manuscript \cite{alpay2003spectral}). Since for the unitary shift,
  $\mathrm Z^{*m}\mathcal Z^n=\mathcal Z^{n-m}$,
we have
\[
\frac{1}{2\pi i}\int_{|z|=r}D_m^*\mathcal Z^{*m}\overline{z^m}z^n\mathcal Z^nD_n\frac{dz}{z}=\begin{cases}\, 0\quad\quad\hspace{0.9cm}{\rm if}\,\,\, n\not=m,\\
r^{2n}D_n^*D_n\quad{\rm if}\,\,\, n\not=m\end{cases}
\]
we have, using the same methods as above, and with the same notations,
\begin{equation}
  \int_{|z|<R}(U(z))^*U(z)d\mu(z)=\sum_{n=0}^\infty D_n^*D_n.
  \label{inner-prod-matrix-ltv}
\end{equation}
  \end{remark}
\subsection{The Fock space}
The classical Fock space corresponds to $\gamma_n=n!$ in  \eqref{gamma-n}.
See \cite{MR0157250}. Besides Bargmann's celebrated characterization ($M_z^*=\partial$, where
$M_z$ is multiplication by $z$ and $\partial$ denotes differentiation)), the Fock space can be seen as the only Hilbert space of functions analytic in a (say convex open)
neighborhood of the origin in which the backward shift operator
 \begin{equation}
      \label{r00}
      (R_0f)(z)=\begin{cases}\, \dfrac{f(z)-f(0)}{z},\quad z\in V\setminus\left\{0\right\}\\
      \,\,\, f^\prime(0),\quad\hspace{11.3mm} z=0,
      \end{cases}
      \end{equation}
      is bounded and has for adjoint the integration operator
      \begin{equation}
(\mathbf I f)(z)=\int_{[0,z]}f(s)ds.
        \label{inte}
        \end{equation}
See \cite{Trevor23}.
The above motivates the following definitions, which are the counterpart of \eqref{r00} and \eqref{inte} in our present setting.
    \begin{definition}
      Let $F(Z)=\sum_{n=0}^\infty Z^nF_n$ be a matrix power series converging in a neighborhood of $0_{p\times p}$. We define
      \begin{equation}
(R_0F)(Z)=\sum_{n=1}^\infty Z^{n-1}F_n.
\label{R-0}
\end{equation}
\begin{equation}
  \mathbf I\left(\sum_{n=0}^\infty Z^nF_n\right)=\sum_{n=0}^\infty \frac{1}{n+1}Z^nF_n.
  \end{equation}
\end{definition}

It is not difficult to check that
\begin{equation}
  \mathbf I(F)=\int_{[0,z]}F(sZ)ds.
  \end{equation}
    
    \begin{definition}
      The Fock space $\mathfrak F$ consists of the matrix power series $F(Z)=\sum_{n=0}^\infty Z^nF_n$ for which
  \begin{equation}
    {\rm Tr}\,[F,F]_{\mathfrak F}<\infty,
    \end{equation}
    where
    \begin{equation}
[F,F]_{\mathfrak F}=    \sum_{n=0}^\infty n!F_n^*F_n.
\end{equation}
  \end{definition}
    We can apply Theorem \ref{thmeasure} with $m(r)=e^{-r^2}$ and obtain:
    \begin{equation}
\sum_{n=0}^\infty n!F_n^*F_n=\frac{1}{\pi}\iint_{\mathbb C}F(zI_p)^*F(zI_p)e^{-|z|^2}dxdy
\label{new567}
      \end{equation}

For the case $p=1$ the following result has been proved in \cite[Lemma 2.1]{Trevor23}.
    
    \begin{proposition}
      The Fock space functions are defined for all $Z\in\mathbb C^{p\times p}$ and the Fock space can be characterized, up to a positive multiplicative
      factor for the inner product,  as the only space of matrix power series {\it a priori} defined near the origin and for which
        \begin{equation}
          R_0^*=\mathbf I.
          \end{equation}
        \end{proposition}

    \begin{proof}
      The first claim follows from $\|Z^n\|\le \|Z\|^n$. For the second claim, and with the understanding that $R_0C=0$ we have
\[
  \begin{split}
    [R_0^*(Z^nC),Z^mD]_{\mathfrak F}&=\begin{cases}\,\,[Z^nC, Z^{m-1}D]_{\mathfrak F},\,\, m=1,2,\ldots\\
    \,\,  0,\hspace{2.66cm}m=0,\end{cases}\\
    &=\delta_{n,m-1} (m-1)!D^*C\\
    &=\delta_{n+1,m}[\frac{1}{n+1}Z^{n+1}C,Z^mD]_{\mathfrak F}\\
    &=[\mathbf I(Z^nC),Z^mD]_{\mathfrak F}.
          \end{split}
          \]
    \end{proof}

    \subsection{Wiener algebra}
    The Wiener algebra of the unit circle, denoted here $\mathcal W$, consists of the trigonometric series of the form
    \begin{equation}
f(e^{it})=\sum_{n\in\mathbb Z}f_ne^{int}
      \end{equation}
    where the complex numbers $f_n$ satisfy $\sum_{n\in\mathbb Z}|f_n|<\infty$. The celebrated Wiener-L\'evy theorem states that $f$ is invertible in $\mathcal W$
    if and only if it is pointwise invertible (that is, pointwise different from $0$). The Wiener algebra is an example of Banach algebra. Two important subalgebras of the Wiener algebra consists of $\mathcal W_+$ (resp.
    $\mathcal W_-$) for which $f_n=0$, $n<0$ (resp. $f_n=0$, $n>0$). The Wiener algebra still makes sense when the $f_n\in\mathbb C^{p\times p}$.
    The counterpart of the Wiener-L\'evy theorem involves then the determinant of the function.\smallskip

We now define the counterpart of $\mathcal W_+$  in the present framework.

\begin{definition}
    We denote by     $\mathfrak W_+$ the space of functions
    \[
F(Z)=    \sum_{n\in\mathbb N_0} Z^n F_n
    \]
    with $\rho(Z)\le 1$ and $\sum_{n=0}^\infty \|F_n\|<\infty$.
  \end{definition}

We leave to the reader to check that $\mathfrak W_+$ endowed with the $\star$ product is an algebra.
  
    \begin{theorem}
      $F\in\mathfrak W_+$ is invertible in $\mathfrak W_+$ if and only if
      \[
\det(\sum_{n=0}^\infty z^nF_n)\not=0,\quad |z|\le1.
      \]
    \end{theorem}

    \begin{proof}
      We first note that $F(zI_p)\in\mathcal W^{p\times p}_+$. So it will be invertible in $\mathcal W^{p\times p}_+$ if and only if $\det F(zI_p)\not=0$ in the closed unit disk. Assume this condition in force, and let $g(z)=\sum_{n=0}^\infty z^nG_n\in\mathcal W_+^{p\times p}$ (the lower case $g$ is not a misprint)
      be such that
      \[
\sum_{n=0}^\infty \|G_n\|<\infty
\]
and $F(zI_p)g(z)=I_p$, for $|z|\le 1$. Then, we have
\[
  \begin{split}
    F_0G_0&=I_p\\
    F_0G_1+F_1G_0&=0\\
    &\hspace{2.5mm}\vdots\\
    F_0G_n+\cdots +F_nG_0&=0\\
&    \hspace{2.5mm}\vdots\\
    \end{split}
  \]
  These equalities express that $g(z)=G(zI_p)$ where $G=\sum_{n=0}^\infty Z^nG_n\in \mathfrak W_+$. The converse statement is proved by reading backwards these arguments.
\end{proof}

The definition of the counterpart of $\mathcal W$ will involve a unitary variable $Z$. The case of the counterpart of $\mathcal W_-$ will be more problematic since one then requires invertible $Z$. These aspects will be considered in a different publication, where
a counterpart of the Wiener-L\'evy inversion theorem is also considered.

        \subsection{Rational functions}
    We consider $(\mathbb C^{p\times p})^{u\times v}$ as a right module over $\mathbb C^{p\times p}$ and a left module over $\mathbb C^{p\times p}$ in the following ways:
    \begin{definition}
      Right module structure: For every $u,v\in \mathbb N$, every $A_{jk}\in\mathbb C^{p\times p}$, $j=1,\ldots, u$, $k=1,\ldots, v$, and every $C\in\mathbb C^{p\times p}$
      \begin{equation}
\left((A_{jk})_{\substack{j=1,\ldots u\\ k=1,\ldots v}}\right)C=(A_{jk}C)_{\substack{j=1,\ldots u\\ k=1,\ldots v}}
\end{equation}
Left module structure: in the above notation
\begin{equation}
  C\left((A_{jk})_{\substack{j=1,\ldots u\\ k=1,\ldots v}}\right)=(CA_{jk})_{\substack{j=1,\ldots u\\ k=1,\ldots v}}.
  \end{equation}
      \end{definition}

      We consider power series of the form
      \begin{equation}
\sum_{n=0}^\infty Z^nA_n,\quad A_n\in (\mathbb C^{p\times p})^{u\times v}.
\label{label7890}
      \end{equation}
We extend $R_0$ on these power series by
\[
R_0(\sum_{n=0}^\infty Z^nA_n)=\sum_{n=1}^\infty Z^{n-1}A_n.
  \]
      \begin{proposition}
        Let $\mathfrak M$ be a finitely generated right module of power series of the form \eqref{label7890}. Assume that $\mathfrak M$ is $R_0$-invariant.
Then        $\mathfrak M$ is generated by the block columns of a matrix function of the form
        \[
E(Z)=C\star(I_{pN}-ZA)^{-\star}
\]
where $A\in(\mathbb C^{p\times p})^{N\times N}$ and $C\in(\mathbb C^{p\times p})^{u\times N}$.
\end{proposition}

\begin{proof}
  Let $F_1,\ldots, F_M$ a generating family for $\mathfrak M$. Every $F_j$ is $(\mathbb C^{p\times p})^{u\times v}$-valued and any $F\in\mathfrak M$ can be written (possibly in a non-unique way) as
  \[
F=\sum_{m=1}^MF_mC_m,\quad C_1,\ldots, C_M\in\mathbb C^{p\times p}.
    \]
    Let
    \[
 R_0F_j=\sum_{r=1}^MF_rA_{rj}.
    \]
    Let $E=\begin{pmatrix}F_1&F_2&\cdots &F_M\end{pmatrix}$. Then,
    \[
      \begin{split}
        R_0F&=\begin{pmatrix}R_0F_1&R_0F_2&\cdots &R_0F_M\end{pmatrix}\\
        &=\begin{pmatrix}F_1A_{11}+F_2A_{21}+\cdots          &F_1A_{12}+F_2A_{22}+\cdots&\cdots &F_1A_{1M}+F_2A_{2M}+\cdots\end{pmatrix}\\
        &=\begin{pmatrix}F_1&F_2&\cdots &F_M\end{pmatrix}\underbrace{\begin{pmatrix}A_{11}&A_{12}&\cdots &A_{1M}\\
          A_{21}&A_{22}&\cdots &A_{2M}\\
          & & &\\
          A_{M1}&A_{M2}& &A_{MM}
        \end{pmatrix}}_{A\in(\mathbb C^{p\times p})^{M\times m}}.
      \end{split}
    \]
    Let $E=\sum_{n=0}^\infty Z^nE_n,\quad E_n\in(\mathbb C^{p\times p})^{M\times M}$. The above equation can be rewritten as
    \[
      \sum_{n=1}^\infty Z^{n-1}E_n=\sum_{n=0}^\infty Z^nE_nA.
    \]
    Hence
    \[
      E_{n+1}=E_nA,\quad n=0,1,\ldots
    \]
    and so $E_n=E_0A^n$, $n=0,1,\ldots$.
    Thus
    \[
      E(Z)=\sum_{n=0}^\infty Z^nE_0A^n=E_0\star(I_{Np}-ZA)^{-\star}.
    \]
\end{proof}
    \begin{definition}
$E$ of the form \eqref{label7890} is rational if the right linear module over $\mathbb C^{p\times p}$ generated by $R_0^jE$, $j=1,\ldots$ is finitely generated.
      \end{definition}

      \begin{theorem}
        $E$ is rational if and only if
        \begin{equation}
          \label{M-real}
          E(Z)=D+Z C\star(I_{pN}-ZA)^{-\star}B.
        \end{equation}
      \end{theorem}
      \begin{proof}
       Let  $\mathfrak M$ be the module generated by $R_0^jE$. Then $\mathfrak M$ is generated by $C\star(I_{pN}-ZA)^{-\star}$. Write
        \[
          R_0E=C\star(I_{PN}-ZA)^{-\star}B.
        \]
        So
        \[
ZR_0E=E-E_0=Z C\star(I_{PN}-ZA)^{-\star}B,
\]
and hence the result holds.
\end{proof}

We note that $E(zI_p)$ is rational in the classical sense and the restriction of \eqref{M-real} to $Z=zI_p$ gives a realization in the classical sense; see
\cite{MR2363355}.

\section{The Hardy space}
\setcounter{equation}{0}
\subsection{Definition}
The classical Hardy space $\mathbf H_2(\mathbb D)$ of the open unit disk $\mathbb D$ is the space of power series $f(z)=\sum_{n=0}^\infty z^n f_n$ for which
\[
  \sum_{n=0}^\infty|f_n|^2<\infty,
  \]
i.e., corresponding to $\gamma_n\equiv 1$ in \eqref{gamma-n}. It is the reproducing kernel Hilbert space with reproducing kernel
\begin{equation}
  \label{szego}
  k(z,w)=\frac{1}{1-z\overline{w}}=\sum_{n=0}^\infty z^n\overline{w}^n,
  \end{equation}
where $z,w$ run through $\mathbb D$, and plays a key role in operator theory. It has numerous extensions and generalizations. In the present paper we consider its counterpart when the complex numbers are replaced by elements in $\mathbb C^{p\times p}$, and 
the Szeg\"o kernel \eqref{szego} replaced by 
\begin{equation}
  \label{matrix-szego}
K(Z,W)=\sum_{n=0}^\infty Z^nW^{*n}
  \end{equation}
  where $Z,W$ runs through the set of elements of $\mathbb C^{p\times p}$ with spectral radius less than $1$.\\

\begin{definition}
   We denote by
\begin{equation}
  \mathbb K=\left\{Z\in\mathbb C^{p\times p}\,;\,\rho(Z)<1\right\}
\end{equation}
where $\rho(Z)$ denotes the spectral radius of $Z$.
\end{definition}
\begin{definition}
  The Hardy space $\mathbf H_2(\mathbb K)$  consists of functions of the form
\begin{equation}
  F(Z)=\sum_{n=0}^\infty Z^nF_n,\quad Z\in\mathbb K,
\label{fz}
\end{equation}
where $(F_0,F_1,\ldots)$ satisfy
\begin{equation}
  \label{normnorm}
 {\rm Tr}\,\left(\sum_{n=0}^\infty F_n^*F_n\right)<\infty.
\end{equation}
\end{definition}

\begin{theorem}
When endowed with the $\mathbb C^{p\times p}$-valued Hermitian form
  \begin{equation}
[F,G]_2=\sum_{n=0}^\infty G_n^*F_n,
\label{HF}
  \end{equation}
(where $G(Z)=\sum_{n=0}^\infty Z^nG_n$)
and associated norm \eqref{normnorm}, $\mathbf H_2(\mathbb K)$  the reproducing kernel Hilbert module with reproducing kernel
\eqref{matrix-szego}, meaning that
\begin{equation}
[F(\cdot), K(\cdot, W)C]=C^*F(W),\quad W\in\mathbb K\,\,\,{\rm and}\,\,\, C\in\mathbb C^{p\times p}.
\label{RK}
  \end{equation}
\end{theorem}

\begin{proof}

We consider the Hilbert space $\ell_2(\mathbb N_0,\mathbb C^{p\times p})$ of sequences $\mathbf F=(F_0,F_1,\ldots)$ of sequences of elements
of $\mathbb C^{p\times p}$ with finite norm \eqref{normnorm}. Note that $\ell_2(\mathbb N_0,\mathbb C^{p\times p})$ is a
right $\mathbb C^{p\times p}$ module: it holds that
\[
\mathbf  F=(F_0,F_1,\ldots)\in\ell_2(\mathbb N_0,\mathbb C^{p\times p})\quad\Longrightarrow\quad \mathbf F A=(F_0A,F_1A,\ldots)
\in\ell_2(\mathbb N_0,\mathbb C^{p\times p}),\quad\forall A\in\mathbb C^{p\times p},
\]
as well as the linearity conditions
\[
\begin{split}
  (\mathbf F+\mathbf G)A&=\mathbf F A+\mathbf G A\\
  \mathbf F (A+B)&=  \mathbf F A+  \mathbf F B,
\end{split}
\]
for all $\mathbf F,\mathbf G\in\ell_2(\mathbb N_0,\mathbb C^{p\times p})$ and $A,B\in\mathbb C^{p\times p}$; see e.g. \cite[p. 117]{godement}.\smallskip

We define
\begin{equation}
  [\mathbf C,\mathbf D]=\sum_{n=0}^\infty G_n^*F_n,\quad \mathbf C\,\,\,{\rm and}\,\,\,{\mathbf D}\in \ell_2(\mathbb N_0,\mathbb C^{p\times p}).
\end{equation}
It holds that
\begin{equation}
  {\rm Tr}\,[\mathbf C,\mathbf D]=\sum_{n=0}^\infty {\rm Tr}\,G_n^*F_n,
  \end{equation}
which is the inner product associated to \eqref{normnorm}.\smallskip

  $\mathbf H_2(\mathbb K)$ is a Hilbert space since $\ell_2(\mathbb N_0,\mathbb C^{p\times p})$ is a Hilbert space, and since the coeffiicents $F_n$ in \eqref{fz} are uniquely determined by $F(Z)$ (and in fact by $F(z I_p)$ with $z\in\mathbb C$).
  Furthermore we have with $F$ as in \eqref{fz},
  \[
    \begin{split}
      [F(\cdot), K(\cdot, W)C]_{\mathbf H_2(\mathbb K)}      &=[\sum_{n=0}^\infty Z^nF_n,\sum_{n=0}^\infty Z^nW^{*n}C]_{\mathbf H_2(\mathbb K)}\\
      &=\sum_{n=0}^\infty C^*W^nF_n\\
      &=[F(W),C]_{\mathbb C^{p\times p}}.
      \end{split}
      \]
\end{proof}

\begin{theorem}
  \label{theorem-inner-product}
  Let $F=\sum_{n=0}^\infty Z^nF_n\in\mathcal H_1$. It holds that
  \begin{equation}
    \lim_{\substack{r\uparrow 1\\r\in(0,1)}}\frac{1}{2\pi}\int_0^{2\pi}F(re^{i\theta}I_p)^*F(re^{i\theta}I_p)d\theta=\sum_{n=0}^\infty F_n^*F_n
  \end{equation}
  where both sides simultaneously converge or diverge.
\end{theorem}

\begin{proof} 
The arguments as in the proof of  Theorem \ref{thmeasure} show that
  \[
\frac{1}{2\pi}\int_0^{2\pi}F(re^{i\theta}I_p)^*F(re^{i\theta}I_p)d\theta=\sum_{n=0}^\infty r^{2n}F_n^*F_n.
\]
Thus, for every $w\in\mathbb C^p$,
  \[
\frac{1}{2\pi}\int_0^{2\pi}w^*F(re^{i\theta}I_p)^*F(re^{i\theta}I_p)w d\theta=\sum_{n=0}^\infty r^{2n}w^*F_n^*F_nw.
\]
Then use the monotone convergence theorem and the polarization identity
\begin{equation}
\label{polar456}
w^*Xz=\frac{1}{4}\sum_{k=0}^3i^{-k}(w+i^kz)^*X(w+i^kz),\quad w,z\in\mathbb C^p,\quad X\in\mathbb C^{p\times p}.
\end{equation}
  \end{proof}

In what follows we will make frequent use of the operators
  \begin{eqnarray}
    (M_ZF)(Z)&=&\sum_{n=0}^\infty Z^{n+1}F_n,\\
    (M_AF)(Z)&=&\sum_{n=0}^\infty Z^nAF_n.
                 \label{maleft}
    \end{eqnarray}
Note that $M_A$ satisfies
\begin{equation}
      \label{MA*}
      (M_A)^*=M_{A^*}. 
      \end{equation}
      Note also that $M_Z$ is an isometry with adjoint $R_0$:
        \begin{equation}
      (M_Z^*)\left(\sum_{n=0}^\infty Z^nF_n\right)=\sum_{n=1}^\infty Z^{n-1}F_n.
      \label{R0}
\end{equation}
      

    
  

    Note that the right hand side of \eqref{R0} makes sense for converging matrix power series even if the series does not belong to
    the Hardy space. Furthermore, for $p=1$ the formula reduces to the classical backward-shift operator defined by \eqref{r00} for a function $f$ analytic 
    in a neighborhood $V$ of the origin.

    \subsection{Resolvent operators and resolvent equations}
    We follow \cite[\S 2.3]{MR1776956} suitably adapted to the present setting. Recall that $M_A$ denotes the operator of multiplication of the coefficients
    on the left by $A$; see \eqref{maleft}.
    We define an operator $R_A$ on power series in $Z$ via:
    \begin{equation}
      (R_AF)(zI_p)=(zI_p-A)^{-1}(F(zI_p)-F(A)).
      \label{R-A}
    \end{equation}
When $A=0_{p\times p}$ this is the operator $R_0$ defined in \eqref{R-0}.
    \begin{lemma}
      Let $F=\sum_{n=0}^\infty Z^nF_n$. Then,
      \begin{equation}
R_AF=\sum_{k=0}^\infty Z^k\left(\sum_{n=k+1}^\infty A^{n-1-k}F_n\right)
        \end{equation}
    \end{lemma}

    \begin{proof}
    We have
    \[
\begin{split}
  R_AF(zI_p)&=(zI_n-A)^{-1}\left(\sum_{n=1}^\infty (z^nI_p-A^n)F_n\right)\\
    &=\sum_{n=0}^\infty (zI_p-A)^{-1}(z^nI_p-A^n)F_n\\
    &=\sum_{n=0}^\infty \sum_{k=0}^{n-1-k}(z^kA^{n-1-k})F_n\\
    &=\sum_{k=0}^\infty z^k\left(\sum_{n=k+1}^\infty A^{n-1-k}F_n\right)
  \end{split}
\]
and hence the result holds.
\end{proof}

\begin{lemma} (see also \cite[(2.32) p. 265]{MR1776956} for the time-varying counterpart) 
\begin{equation}
R_A=R_0(I-M_AR_0)^{-1}.
  \end{equation}
  \end{lemma}
\begin{proof}
  We follow \cite{MR1776956}  and set $G=(I-M_AR_0)F$. We have:
\[
 \begin{split}
   G&=F-M_AR_0F\\
   &=\sum_{n=0}^\infty Z^nF_n-M_A\sum_{n=1}^\infty Z^{n-1}F_n\\
   &=F_0+\sum_{n=1}^\infty(Z^n-Z^{n-1}A)F_n\\
   &=F_0+\sum_{n=0}^\infty Z^{n-1}(Z-A)F_n
    \end{split}
  \]
  and so $G(A)=F_0$ and
  \begin{equation}
    G(zI_p)=F_0+(zI_p-A)(R_0F)(zI_p).
  \end{equation}
  Thus
  \[
    \begin{split}
      (R_AG)(zI_p)&=(zI_p-A)^{-1}(G(zI_p)-G(A))\\
      &=(zI_n-A)^{-1}((G-F_0)(zI_p))\\
      &=(R_0F)(zI_p)
      \end{split}
    \]
    so that
    \[
      R_A(I-M_AR_0)=R_0
    \]
    and hence the result holds.
  \end{proof}

  We now have the resolvent equation:
  
  \begin{theorem}  Let $A,B\in\mathbb C^{p\times p}$. It holds that:
    \begin{equation}
R_A-R_B=R_A(M_A-M_B)R_B.
      \end{equation}
    \end{theorem}

    \begin{proof}
      Using
      \[
        R_0(I-M_AR_0)^{-1}=(I-R_0M_A)^{-1}R_0
      \]
      we can write
      \[
        \begin{split}
          R_A-R_B&=R_0(I-M_AR_0)^{-1}-R_0(I-M_BR_0)^{-1}\\
          &=(I-R_0M_A)^{-1}R_0-R_0(I-M_BR_0)^{-1}\\
          &=(I-R_0M_A)^{-1}\left(R_0-R_0M_BR_0-R_0+R_0M_AR_0\right)R_0(I-M_BR_0)^{-1}\\
            &=(I-R_0M_A)^{-1}\left(R_0M_AR_0-R_0M_BR_0\right)R_0(I-M_BR_0)^{-1}\\
              &=(I-R_0M_A)^{-1}R_0(M_A-M_B)R_0(I-M_BR_0)^{-1}\\
                  &=R_0(I-M_AR_0)^{-1}(M_A-M_B)R_0(I-M_BR_0)^{-1}\\
        \end{split}
      \]
      and hence the result holds.
      \end{proof}
\subsection{Blaschke factor}
Given a matrix $A\in\mathbb C^{p\times p}$ with $\rho(A)<1$, the unique solution to the Stein equation
\begin{equation} 
 \label{gamma-eq}
  \Gamma_A-A\Gamma_AA^*=I_p
\end{equation}
is given by the converging series 
\begin{equation}
  \Gamma_A=\sum_{n=0}^\infty A^nA^{*n}.
\label{PS}
\end{equation}
Observe that $\Gamma_A\geq I_p$ and hence is invertible. If we let
\[
L_A=\Gamma_A-\Gamma_AA^*\Gamma_A^{-1}A\Gamma_A,
  \]
then it follows by the Sherman-Morrison formula that
\begin{equation}
  \label{eq00005}
L_A^{-1}=A^*A+\Gamma_A^{-1}\ge  0.
    \end{equation}
Let us now introduce the power series 
\begin{align}
U_A(Z)&:=(Z-A)\star (I-Z\Gamma_A A^*\Gamma_A^{-1})^{-\star}L_A^{\frac{1}{2}}\notag\\
&=-AL_A^{\frac{1}{2}}+\sum_{n=1}^\infty Z^n(I-A\Gamma_AA^*\Gamma_A^{-1})(\Gamma_A A^*\Gamma_A^{-1})^{n-1} 
L_A^{\frac{1}{2}}\notag\\
&=-AL_A^{\frac{1}{2}}+\sum_{n=1}^\infty Z^n A^{*(n-1)}\Gamma_A^{-1}L_A^{\frac{1}{2}}.
\label{UA}
\end{align}
The first power series representation of $U_A$ above follows since 
$$
(I-Z\Gamma_A A^*\Gamma_A^{-1})^{-\star}=\sum_{k=0}^\infty Z^k(\Gamma_AA^*\Gamma_A^{-1})^k,
$$
while the next representation follows since 
\begin{align*}
(I-A\Gamma_AA^*\Gamma_A^{-1})(\Gamma_A A^*\Gamma_A^{-1})^{n-1}&=(I-A\Gamma_AA^*\Gamma_A^{-1})\Gamma_A A^{*(n-1)}\Gamma_A^{-1}\\
&=(\Gamma_A-A\Gamma_AA^*)A^{*(n-1)}\Gamma_A^{-1}=A^{*(n-1)}\Gamma_A^{-1}
\end{align*}
for all $n\ge 1$, due to \eqref{gamma-eq}.
\begin{proposition}
Let $U_A(Z)$ be defined as in \eqref{UA}. Then
\begin{equation}
\left[M_Z^n U_A, \, M_Z^k U_A\right]_2=\delta_{n,k}I_p\quad \mbox{for all}\quad n,k\ge 0.
\label{UAnorm}
\end{equation}
\label{P:bf}
\end{proposition}
\begin{proof}
By \eqref{R0}, it suffices to verify \eqref{UAnorm} for $k=0$. For $n=0$, we have by \eqref{UA} and definitions
\eqref{HF} and \eqref{PS},
\begin{align*}
\left[U_A, \, U_A\right]_2&=L_A^{\frac{1}{2}}A^*AL_A^{\frac{1}{2}}+\sum_{j=0}^{\infty}L_A^{\frac{1}{2}}\Gamma_A^{-1}A^{j}A^{*j}\Gamma_A^{-1}
L_A^{\frac{1}{2}}\\
&=L_A^{\frac{1}{2}}\big(A^*A+\Gamma_A^{-1}\bigg(\sum_{j=0}^{\infty}A^{j}A^{*j}\bigg)\Gamma_A^{-1}\big)L_A^{\frac{1}{2}}\\
&=L_A^{\frac{1}{2}}\big(A^*A+\Gamma_A^{-1}\big)L_A^{\frac{1}{2}}=L_A^{\frac{1}{2}}L_AL_A^{\frac{1}{2}}=I_p.
\end{align*}
For $n>0$, we have
$$
Z^nU_A(Z)=-Z^nAL_A^{\frac{1}{2}}+\sum_{j=1}^\infty Z^{j+n} A^{*(j-1)}\Gamma_A^{-1}L_A^{\frac{1}{2}},
$$
and subsequently,
\begin{align*}
\left[M_Z^nU_A, \, U_A\right]_2&=-L_A^{\frac{1}{2}}\Gamma_A^{-1}A^{n-1}AL_A^{\frac{1}{2}}+\sum_{j=0}^{\infty}L_A^{\frac{1}{2}}\Gamma_A^{-1}A^{n+j}A^{*j}
\Gamma_A^{-1}L_A^{\frac{1}{2}}\\
&=L_A^{\frac{1}{2}}\Gamma_A^{-1}A^n\big(-I_p+\bigg(\sum_{j=0}^{\infty}A^{j}A^{*j}\bigg)\Gamma_A^{-1}\big)L_A^{\frac{1}{2}}\\
&=L_A^{\frac{1}{2}}\Gamma_A^{-1}A^n\big(-I_p+\Gamma_A\Gamma_A^{-1}\big)L_A^{\frac{1}{2}}=0,
\end{align*}
which completes the proof.
\end{proof}
\begin{corollary}
The operator $M_{U_A}: \, F(z)\mapsto U_A(z)\star F(z)$ is an isometry on ${\bf H}_2(\mathbb K)$ in the following sense:
\begin{equation}
\left[U_A\star F, \, U_A\star F\right]_2=\left[F, \, F\right]_2\quad\mbox{for any}\quad F\in{\bf H}_2(\mathbb K).
\label{isom1}
\end{equation}
\label{C:ISO}
\end{corollary}
Indeed, if we take $F\in{\bf H}_2(\mathbb K)$ in the form \eqref{fz}, then we have, by \eqref{UAnorm}, 
\begin{align*}
\left[U_A\star F, \, U_A\star F\right]_2&=
\bigg[\sum_{j=0}^\infty M_Z^j U_A C_j, \, \sum_{j=0}^\infty M_Z^j U_A C_j\bigg]_2\\
&=\sum_{j=0}^\infty \left[U_A C_j, \; U_A C_j\right]_2=\sum_{j=0}^\infty \left[C_j, \; C_j\right]_2=
\sum_{j=0}^\infty C_j^*C_j=\left[F, \, F\right]_2.
\end{align*}
\begin{proposition}
  \label{bastille}
$\mathbf H_2(\mathbb K)\ominus M_{U_A}\mathbf H_2(\mathbb K)=\left\{K(Z,A)C:\; C\in\mathbb C^{p\times p}\right\}$.
  \end{proposition}
  \begin{proof}
Since $U_A(A)=0$, it follows by Remark \ref{R:eval} that $(U_A\star F)(A)=0$ for any $F\in \mathbf H_2(\mathbb K)$.
Therefore, by \eqref{RK} we have 
$$
\left[U_A\star F, \, K(Z,A)C\right]_2=C^*(U_A\star F)(A)=0
$$
for all $F\in \mathbf H_2(\mathbb K)$ and $C\in\mathbb C^{p\times p}$. Therefore, the submodules
 $\left\{K(Z,A)C:\; C\in\mathbb C^{p\times p}\right\}$ and 
$M_{U_A}\mathbf H_2(\mathbb K)$ of $\mathbf H_2(\mathbb K)$ are orthogonal with respect to the form \eqref{HF}.

\smallskip

We next take an arbitrary $G(Z)=\sum_{j=0}^\infty Z^jG_j\in \mathbf H_2(\mathbb K)$, which is orthogonal to
$M_{U_A}\mathbf H_2(\mathbb K)$. In particular, it is orthogonal to $Z^nU_A$ for all $n\ge 0$. In other words,
$$
\left[G, \, M_Z^nU_A\right]_2=0\quad\mbox{for all}\quad n\ge 0,
$$ 
which can be written in terms of coefficients of $G$ and $U_A$ (see \eqref{UA}) as 
$$
-L_A^{\frac{1}{2}}A^*G_n+\sum_{j=1}^\infty L_A^{\frac{1}{2}}\Gamma_A^{-1}A^{j-1}G_{n+j}=0,
$$
or equivalently,  
\begin{equation}
-\Gamma_AA^*G_n+\sum_{j=1}^\infty A^{j-1}G_{n+j}=0
\label{aug23b}
\end{equation}
for all $n\ge 0$. Replacing $n$ by $n+1$ in \eqref{aug23b} and multiplying both parts by $A$ on the left we get
$$
-A\Gamma_AA^*G_{n+1}+\sum_{j=1}^\infty A^{j}G_{n+j+1}=0.
$$
Upon subtracting the latter equality from \eqref{aug23b} we get
\begin{align*}
0&=\Gamma_AA^*G_n-A\Gamma_AA^*G_{n+1}-G_{n+1}\\
&=\Gamma_AA^*G_n-(I+A\Gamma_AA^*)G_{n+1}=\Gamma_AA^*G_n-\Gamma_AG_{n+1}.
\end{align*}
Therefore, $G_{n+1}=A^{*n}G_n$ for all $n\ge 0$ and hence $G_n=A^{*n}G_0$, so that 
$$
G(z)=\sum_{j=0}^\infty Z^jA^{*n}G_0=K(Z,A)G(0),
$$
which completes the proof.
    \end{proof}
\begin{corollary}
 An element $F$ of $\mathbf H_2(\mathbb K)$ satisfies $F(A)=0$ if and only if
  it is of the form $F=U_A\star G$ for some $G\in\mathbf H_2(\mathbb K)$.
\label{C:hom}
\end{corollary}
\begin{proof}
The ``if" part follows from Remark \ref{R:eval} since $U_A(A)=0$.
On the other hand, if $F\in \mathbf H_2(\mathbb K)$ is subject to $F(A)=0$, then
$$
 [F, K(\cdot, A)C]_2=C^*F(A)=0,\quad  C\in\mathbb C^{p\times p},
$$
and hence, $F\in M_{U_A}\mathbf H_2(\mathbb K)$, by Proposition \ref{bastille}.
\end{proof}

The preceding results are of special interest when $Z=zI_p$. The function $b_A(z)\stackrel{\rm def.}{=}U_A(zI_p)$ is then a $\mathbb C^{p\times p}$-valued rational function contractive in the open unit disk and unitary on the unit circle. As such one can apply the results of \cite{ad3,ag} which characterize such
functions, as we now explain. We rewrite first
\begin{eqnarray}
  b_A(z)&=&-AL_A^{1/2}+z(I_p-zA^*)^{-1}\Gamma_A^{-1}L_A^{1/2}\\
        &=&(zI_p-A)(I_p-z\Gamma_AA^*\Gamma_A^{-1})^{-1}L_A^{1/2}\\
            &=&(I_p-zA^*)^{-1}(zI_P-AL_A^{-1})L_A^{1/2}.
  \end{eqnarray}
  These equations are a direct consequence of \eqref{UA} when $Z=zI_p$. When $A=aI_p$ with $a\in\mathbb D$, we have
  \[
b_A(z)=\frac{z-a}{1-z\overline{a}}I_p.
    \]
    Furthermore, we set
      \begin{equation}
        \label{real-bA}
    \begin{split}
      \mathcal A&=A^*\\
      \mathcal B&=\Gamma_A^{-1}L_A^{1/2}\\
      \mathcal C&=I_p\\
      \mathcal D&=-AL_A^{1/2}
      \end{split}
    \end{equation}
    so that
    \[
b_A(z)=\mathcal D+z\mathcal C(I_p-z\mathcal A)^{-1}\mathcal B.
      \]

    \begin{theorem}
      \eqref{real-bA} is a minimal realization of $b_A(z)$ and it holds that
      \begin{equation}
        \label{superbowl}
        \begin{pmatrix}\mathcal A&\mathcal B\\ \mathcal C&\mathcal D\end{pmatrix}^*
\begin{pmatrix}\Gamma_A&0\\0&I_p\end{pmatrix}
\begin{pmatrix}\mathcal A&\mathcal B\\ \mathcal C&\mathcal D\end{pmatrix}=
\begin{pmatrix}\Gamma_A&0\\0&I_p\end{pmatrix}.
\end{equation}
\label{supersuper}
  \end{theorem}

  \begin{proof} 
    We need to verify that
    \begin{eqnarray}
      A\Gamma_AA^*+I_p&=&\Gamma_A\\
      L_A^{1/2}\Gamma_A^{-1}\Gamma_AA^*-L_A^{1/2}A^*&=&0\\
      L_A^{1/2}\Gamma_A^{-1}\Gamma_A\Gamma_A^{-1}L_A^{1/2}+L_A^{1/2}A^*AL_A^{1/2}&=&I_p.
    \end{eqnarray}
    The  first identity is just \eqref{gamma-eq}, the second one is trivial and the last identity follows from \eqref{eq00005}.
  \end{proof}

  \begin{remark}
In the language of \cite{ag}, $\Gamma_A$ is the associated Hermitian matrix associated to the minimal realization \eqref{real-bA}.
\end{remark}

As a corollary of the previous arguments we see that an element $f(z)=\sum_{n=0}^\infty z^nF_n\in\mathbf H_2(\mathbb D)^{p\times p}$, with $F_0,F_1,\ldots\in\mathbb C^{p\times p}$ satisfies $f(A)=0$ if and only if it can be written as
\[
  f(z)=b_A(z)g(z)
\]
where $\mathbf H_2(\mathbb D)^{p\times p}$, and $\|f\|=\|g\|$. One could also obtain this result directly from \cite{MR96h:47020} or \cite{bgr} for instance.\\

{\bf Particular cases}:
%
%
In the case of quaternions, we take 
 \begin{equation}
Z=\begin{pmatrix}z_1&-z_2\\ \overline{z}_2&\overline{z}_1\end{pmatrix}\quad\mbox{and}\quad 
A=\begin{pmatrix}a_1&-a_2\\ \overline{a}_2&\overline{a}_1\end{pmatrix}, \; \; |a_1|^2+|a_2|^2<1.
     \label{eq-quart}
    \end{equation}
Since $AA^*=A^*A=(|a_1|^2+|a_2|^2)I_2$, we derive from \eqref{PS} $\Gamma_a=(1-|a_1|^2-|a_2|^2)^{-1}$, $L_A=I_2$, and then
\eqref{UA} amounts to
\[
  U_A=(Z-A)\star(I_2-ZA^*)^{-\star}.
\]
See \cite{acs1}, where the quaternionic notation rather than matrix notation is used.\smallskip

The split quaternions correspond to
\[
A=\begin{pmatrix}a_1&a_2\\ \overline{a_2}&\overline{a_1}\end{pmatrix},\quad a_1,a_2\in\mathbb C.
\]
See \cite{alss_IJM,MR2822209,MR3026140}.
We now have, with
\[
  J=\begin{pmatrix}1&0\\0&-1\end{pmatrix},
  \]
  \[
    AJA^*=(|a_1|^2-|a_2|^2)J.
\]
Thus, an indefinite metric appears and the current theory has to be appropriately extended.

\subsection{A Beurling-Lax type theorem}
Let $\mathfrak M$ be a closed subspace of $\mathbf H_2(\mathbb K)$ invariant under $M_Z$. Then, the set of functions of a complex variable $F(zI_p)$ with $F\in\mathfrak M$ is a $z$-invariant subspace of the classical Hardy space $\mathbf H_2(\mathbb D)\otimes\mathbb C^{p\times p}$
(i.e. $\mathbb C^{p\times p}$-valued functions with entries in $\mathbf H_2(\mathbb D)$). By the classical Beurling-Lax theorem there is a Hilbert space $\mathcal C$ and a $L(\mathcal C,\mathbb C^{p\times p})$-valued function $\Theta$ which takes coisometric values on the unit circle and such that
\[
\mathfrak M=\Theta\mathbf H_2(\mathbb D,\mathcal C).
\]
Let $\Theta(z)=\sum_{n=0}^\infty  z^n\Theta_n$ and for $H\in\mathbf H_2(\mathbb D,\mathcal C)$, set $H(z)=\sum_{n=0}^\infty z^nH_n$.
Let $F=\Theta H$, with $F(z)=\sum_{n=0}^\infty z^nF_n$.
{\it A priori} $\mathcal C$ is not a left $\mathbb C^{p\times p}$-module and we cannot lift directy to the setting of the paper, i.e., cannot write from
\begin{equation}
  \label{bekholzot}
    F(z)=\sum_{n=0}^\infty z^nF_n=\Theta(z)H(z)
  \end{equation}
  \[
    \sum_{n=0}^\infty Z^nF_n=\left(\sum_{n=0}^\infty Z^n\Theta_n\right)\star\left(\sum_{n=0}^\infty Z^nH_n\right).
  \]
But \eqref{bekholzot} can be rewritten as 
\begin{equation}
\label{bekholzot!}
F(z)=\sum_{n=0}^\infty z^nF_n=\sum_{n=0}^\infty z^n\left(\sum_{k=0}^n\Theta_kH_{n-k}\right)
\end{equation}  
The operator products $\Theta_kH_{n-k}$ make sense by construction, and \eqref{bekholzot!} is equivalent to
\begin{equation}
\label{bekholzot!!}
\sum_{n=0}^\infty Z^nF_n=\sum_{n=0}^\infty Z^n\left(\sum_{k=0}^n\Theta_kH_{n-k}\right).
\end{equation}  
Hence:

\begin{theorem}
  Let $\mathfrak M$ be a closed subspace of $\mathbf H_2(\mathbb K)$ invariant under $M_Z$. There exist a Hilbert space $\mathcal C$ and an
  operator-valued  function $\Theta=\sum_{n=0}^\infty Z^n\Theta_n$ such that $F=\sum_{n=0}^\infty Z^nF_n  \in\mathfrak M$ if and only if
  \eqref{bekholzot!} holds for some $H\in\mathbf H_2(\mathbb D,\mathcal C)$.
\end{theorem}
\section{Interpolation}
\setcounter{equation}{0}
\label{sec-4}
We want to solve:

\begin{problem}
  \label{illustration-pb}
  Given $A_1,B_1,\ldots, A_N,B_N\in\mathbb C^{p\times p}$    (the interpolation data), describe the set of all functions $F\in\mathbf H_2(\mathbb K)$
  such that
  \begin{equation}
    F(A_j)=B_j,\quad j=1,\ldots, N
    \label{inter}
    \end{equation}
\end{problem}

We follow the classical
approach to interpolation in reproducing kernel spaces (and, more generally, modules; see \cite{ab6}). We look for a solution of the form
\[
  F=\sum_{j=1}^N (I-ZA_j^*)^{-\star}C_j
\]
where $C_1,\ldots, C_N\in\mathbb C^{p\times p}$ are to be found.
Since
\[
F=\sum_{j=1}^N\sum_{n=0}^\infty Z^nA_j^{*n}C_j
  \]
the interpolation conditions lead to
  \[
    B_k=F(A_k)  =\sum_{j=1}^N\sum_{n=0}^\infty A_k^nA_j^{*n}C_j,\quad k=1,\ldots, N,
  \]
 and so
\begin{equation}
  \label{matrix-G}
  \underbrace{\begin{pmatrix}\sum_{n=0}^\infty A_1^nA_1^{*n}&\sum_{n=0}^\infty A_2^nA_1^{*n}&\cdots&\sum_{n=0}^\infty A_N^nA_1^{*n}\\
        \sum_{n=0}^\infty A_2^nA_1^{*n}&\sum_{n=0}^\infty A_2^nA_2^{*n}&\cdots&\sum_{n=0}^\infty A_N^nA_2^{*n}\\
            & & & &\\
            \sum_{n=0}^\infty A_N^nA_1^{*n}&\sum_{n=0}^\infty A_N^nA_2^{*n}&\cdots&\sum_{n=0}^\infty A_N^nA_N^{*n}\end{pmatrix}}_{\mathbf G}
                \begin{pmatrix}C_1\\ C_2\\ \vdots \\ C_N\end{pmatrix}=\begin{pmatrix}B_1\\ B_2\\ \vdots \\ B_N\end{pmatrix}
              \end{equation}
The Gram matrix $\mathbf G$ is positive semi-definite. When it is positive definite, one can solve and get the $C_j$. The fact that $\mathbf G>0$ means that the interpolation points are ``far away'' enough one from the other. Similar phenomena occurs in the time-varying setting. See \cite{DD-ot56}.\\

By the Beurling-Lax theorem considered in the previous section one can then consider the functions for which interpolation is met with $B_1=\cdots =B_N=0$ and
get a description of all solutions, which we now present.\smallskip

Let
\begin{equation}
\label{a-c}
\mathcal C=\underbrace{\begin{pmatrix}I_p&I_p&\cdots&I_p\end{pmatrix}}_{N\,\, {\rm times}}\quad{\rm and}\quad \mathcal A={\rm diag}(A_1^*,A_2^*,\ldots, A_N^*).
\end{equation}
The matrix $\mathbf G$ satisfies
the Stein equation
\begin{equation}
  \mathbf G-\mathcal A^*\mathbf G \mathcal A=\mathcal C^*\mathcal C
\label{structure-matrix}
\end{equation}

\begin{proposition}
Define (with $\star$ product performed block-wise)
\begin{equation}
  \begin{split}
    \Theta(Z)&=I_p-
    (I_p-Z)\star \begin{pmatrix}(I_p-ZA_1^*)^{-\star}&\cdots& (I_p-ZA_N^*)^{-\star}\end{pmatrix}
     \mathbf G^{-1}\begin{pmatrix}(I_p-A_1)^{-1}\\
  \vdots\\
  (I_p-A_N)^{-1}    \end{pmatrix}
  \end{split}
\end{equation}
It holds that
\begin{equation}
  \Theta(A_j)=0,\quad j=1,\ldots, N.
  \end{equation}
\end{proposition}

\begin{proof}
  We have
  \[
    \begin{split}
(I_p-Z)\star C\star
  \begin{pmatrix}(I_p-ZA_1^*)^{-\star}&0&0&\cdots &0\\
    0&    (I_p-ZA_2^*)^{-\star}&0&\cdots &0\\
    0& 0&\ddots&&0\\
    0&0&\cdots&&    (I_p-ZA_N^*)^{-\star}    \end{pmatrix}(A_1)&=\\
  \\
  &\hspace{-14.cm}=\underbrace{\begin{pmatrix} I_p-Z&I_p-Z& \cdots\end{pmatrix}}_{N\,\,{\rm times}}\star
  \begin{pmatrix}\sum_{n=0}^\infty Z^nA_1^{n*}&0&0&\cdots &0\\
    0&  \sum_{n=0}^\infty Z^nA_2^{n*}&0&\cdots &0\\
    0& 0&\ddots&&0\\
    0&0&\cdots&&   \sum_{n=0}^\infty Z^nA_N^{n*}    \end{pmatrix}(A_1)\\
  \\
  &\hspace{-14cm}=  \begin{pmatrix}\sum_{n=0}^\infty (Z^n-Z^{n+1})A_1^{n*}&  \sum_{n=0}^\infty (Z^n-Z^{n+1})A_2^{n*}&\cdots
&   \sum_{n=0}^\infty (Z^n-Z^{n+1})A_N^{n*}    \end{pmatrix}(A_1)\\
  \\
  &\hspace{-14cm}=\begin{pmatrix}(I_p-A_1)G_{11}&(I_p-A_1)G_{12}&\cdots&(I_p-A_1)G_{1N}\end{pmatrix}\\ \\
  &\hspace{-14cm}=(I_p-A_1)\begin{pmatrix}G_{11}&G_{12}&\cdots&G_{1N}\end{pmatrix}.
\end{split}
\]
Hence (here too with $\times$ denoting regular matrix multiplication)
\[
  \begin{split}
    \Theta(A_1)&=I_p-(I_p-A_1)\begin{pmatrix}G_{11}&G_{12}&\cdots&G_{1N}\end{pmatrix}\mathbf G^{-1}\times\\
    &\hspace{5mm}\times
  \begin{pmatrix}(I_p-A_1)^{-1}&0&0&\cdots &0\\
    0&    (I_p-A_2)^{-1}&0&\cdots &0\\
    0& 0&\ddots&&0\\
    0&0&\cdots&&    (I_p-A_N)^{-1}    \end{pmatrix}C^*\\ \\
  &=I_p-(I_p-A_1)\begin{pmatrix}I_p&0&\cdots&0\end{pmatrix}\times\\
    &\hspace{5mm}\times
  \begin{pmatrix}(I_p-A_1)^{-1}&0&0&\cdots &0\\
    0&    (I_p-A_2)^{-1}&0&\cdots &0\\
    0& 0&\ddots&&0\\
    0&0&\cdots&&    (I_p-A_N)^{-1}    \end{pmatrix}C^*\\ \\
  &=I_p-\begin{pmatrix}I_p&0&\cdots&0\end{pmatrix}C^*\\
  &=0.
    \end{split}
  \]
  The same argument works of course for $A_2,\ldots, A_N$.
  \end{proof}

  We now set
\[
  \psi(z)=\Theta(zI_p)=I_p-(1-z)C(I_{Np}-zA)^{-1}\mathbf G^{-1}(I_{Np}-A^*)^{-1}C^*
\]

\begin{proposition}
  It holds that
  \begin{equation}
    \label{ran-theta}
\frac{I_p-\psi(z)\psi(w)^*}{1-z\overline{w}}=C(I_{Np}-zA)^{-1}\mathbf G^{-1}(I_{Np}-\overline{w}A^*)^{-1}C^*
    \end{equation}
    and in particular $\psi$ is a rational inner function and so the operation of multiplication by $\psi$ is an isometry from $\mathbf H_2(\mathbb D,\mathbb C^{p\times p})$ into
    itself.
\end{proposition}
\eqref{ran-theta} is a classical computation based on the identity \eqref{structure-matrix},
which originates with L. de Branges' work; see \cite{MR0229011} for the
latter (see also \cite{Dym_CBMS}). We refer to \cite[Exercise 7.1.17 p. 375 and p. 402]{CAPB_2} for a recent presentation of the computation. Rather than proving \eqref{ran-theta} we present a minimal realization of $\psi$ and compute its associated Hermitian matrix, as
in Theorem \ref{supersuper}.

\begin{theorem}
  With $\mathcal A$ and $\mathcal C$ as in \eqref{a-c} and
  \begin{eqnarray}
    \mathcal B&=&\mathbf G^{-1}(I_{Np}-\mathcal A^*)^{-1}\mathcal C^*\\
    \mathcal D&=&I_p-\mathcal C\mathbf G^{-1}(I_{Np}-\mathcal A^*)^{-1}\mathcal C^*
  \end{eqnarray}
  we have
  \[
    \psi(z)=\mathcal D+z\mathcal C(I_{Np}-z\mathcal A)^{-1}\mathcal B
    \]
    and
     \begin{equation}
        \label{superbowl2}
        \begin{pmatrix}\mathcal A&\mathcal B\\ \mathcal C&\mathcal D\end{pmatrix}^*
\begin{pmatrix}\mathbf G&0\\0&I_p\end{pmatrix}
\begin{pmatrix}\mathcal A&\mathcal B\\ \mathcal C&\mathcal D\end{pmatrix}=
\begin{pmatrix}\mathbf G&0\\0&I_p\end{pmatrix}.
\end{equation}
  \end{theorem}

  \begin{proof} We set $A={\rm diag}(A_1,A_2,\ldots, A_N)$. The $(2,1)$ block equality in \eqref{superbowl2} is
    \[
      \mathcal A^*\mathbf G\mathcal A+\mathcal C^*\mathcal C=\mathbf G
    \]
    which is \eqref{structure-matrix}. The $(1,2)$ block amounts to
    \[
      \mathcal A^*\mathbf G\mathcal B+\mathcal C^*\mathcal D=0,
    \]
    which is equivalent to
    \[
      A\mathbf G(I_{Np}-A)^{*}\mathbf G^{-1}(I_{Np}-A)^{-1}\mathcal C+\mathcal C^*\left(I_p-\mathcal C\mathbf G^{-1}(I_{Np}-A)^{-1}\mathcal C^*\right)=0.
    \]
    Equivalently
    \[
\mathcal C^*+\left(A\mathbf G(I_{Np}-A^*)+\mathcal C^*\mathcal C\right)\mathbf G^{-1}(I_{Np}-A)^{-1}\mathcal C^*=0,
\]
i.e, after using \eqref{structure-matrix}
\[
\mathcal C^*+ (A-I_{Np})\mathbf G\mathbf G^{-1}(I_{Np}-A)^{-1}\mathcal C^*=0,
\]
which clearly holds.  To conclude  we verify that identity holds in the $(2,2)$-block, i.e. that we have
\[
  \begin{split}
  I_p=  \mathcal C(I_{Np}-A^*)^{-1}\mathbf G^{-1}(I_{Np}-A)\mathbf G(I_{Np}-A^*)\mathbf G^{-1}(I_{Np}-A)^{-1}\mathcal C^*+&\\
    &\hspace{-8cm}+I_p-\mathcal C(I_{Np}-A^*)^{-1}\mathbf G^{-1}\mathcal C^*-\mathcal C\mathbf G^{-1}(I_{Np}-A)^{-1}\mathcal C^*+\\
    &\hspace{-8cm}+\mathcal C(I_{Np}-A^*)^{-1}\mathbf G^{-1}\mathcal C^*\mathcal C\mathbf G^{-1}(I_{Np}-A)^{-1}\mathcal C^*.
    \end{split}
  \]
  This amounts to check that
  \[
    \mathcal C(I_{Np}-A^*)^{-1}\mathbf G^{-1}\Delta\mathbf G^{-1}(I_{Np}-A)^{-1}\mathcal C^*=0
  \]
  where
  \[
    \Delta=(I_{Np}-A)\mathbf G(I_{Np}-A^*)+\mathcal C^*\mathcal C-(I_{Np}-A)\mathbf G-\mathbf G(I_{Np}-A^*)=0.
  \]
  But it is readily seen that $\Delta=0$ using \eqref{structure-matrix}.
\end{proof}

We now come back to the original interpolation problem.

\begin{proposition}
  Let $A_1,A_2,\ldots, A_N\in \mathbb K$ be such that the matrix $\mathbf G$ (defined by \eqref{matrix-G}) is strictly positive. Then,
$F\in\mathbf H_2(\mathbb K)$ vanishes at $A_1,A_2,\ldots, A_N$ if and only if it is in the range of $M_\Theta$.
\end{proposition}
\begin{proof}

    It follows by the characterization \eqref{theorem-inner-product} of the inner product in $\mathbf H_2(\mathbb K)$ that the operator $M_\Theta$ of star-multiplication on the left by $\Theta$
    is an isometry from $\mathbf H_2(\mathbb K)$ into itself.
    Thus (and taking into account that $M_\Theta$ is an isometry)
    \[
      \begin{split}
        \mathbf H_2(\mathbb K)&={\rm Ran}\,(I-M_\Theta M_\Theta^*)\oplus{\rm Ran}\,(M_\Theta M_\Theta^*)\\
        &={\rm Ran}\,(I-M_\Theta M_\Theta^*)\oplus{\rm Ran}\,M_\Theta.
      \end{split}
    \]

    To characterize ${\rm Ran}\,(I-M_\Theta M_\Theta^*)$ we first remark that
    $\frac{I_p-\psi(z)\psi(w)^*}{1-z\overline{w}}$ is the reproducing kernel of $\mathbf H_2(\mathbb D,\mathbb C^{p\times p})
    \ominus    \psi\mathbf H_2(\mathbb D,\mathbb C^{p\times p})$ and that
    \begin{equation}
      \mathbf H_2(\mathbb D,\mathbb C^{p\times p})    \ominus    \psi\mathbf H_2(\mathbb D,\mathbb C^{p\times p})={\rm Ran}\,(I-M_\psi M_\psi^*).
      \end{equation}
      From \eqref{ran-theta} follows that ${\rm Ran}\,(I-M_\psi M_\psi^*)$ is spanned by the functions $(I_p-zA_j^*)^{-1}$ with coefficients on the right belonging to
      $\mathbb C^{p\times p}$, and this ends the proof since
      \begin{equation}
        \label{functions-!}
       \left( (I_p-ZA_j^*)^{-\star}D\right)(zI_p)=(I_p-zA_j^*)^{-1}D, \quad j=1,\ldots, N\quad{\rm and}\quad D\in\mathbb C^{p\times p}.
\end{equation}
This concludes the proof since a function $F\in\mathbf H_2(\mathbb K)$ vanishes at $A_1,\ldots, A_N$ if and only if it is orthogonal to the functions \eqref{functions-!}.
\end{proof}

Combining with the beginning of the section we have:

\begin{theorem}
  Let $A_1,A_2,\ldots, A_N\in \mathbb K$ be such that the matrix $G$ (defined by \eqref{matrix-G}) is strictly positive, and let $B_1,\ldots, B_N\in\mathbb C^{p\times p}$.
  Then $F$ satisfies the interpolation conditions \eqref{inter}
  \begin{equation}
    F(A_j)=B_j,\quad j=1,\ldots, N
    \label{inter*}
    \end{equation}
if and only F is of the form
\begin{equation}
  \label{ortho00}
F=F_{\rm min}+\Theta\star G
\end{equation}
where
\[
  F_{\rm min}=\begin{pmatrix}I_p-ZA_1^*)^{-\star}&(I_p-ZA_2^*)^{-\star}&\cdots&(I_p-ZA_N^*)^{-\star}\end{pmatrix}
  \mathbf G^{-1}\begin{pmatrix}B_1\\ \vdots \\ B_N\end{pmatrix}
\]
and $G$ runs through $\mathbf H_2(\mathbb K)$. The decomposition \eqref{ortho00} is orthogonal.
\end{theorem}

\section{Schur multipliers}
\setcounter{equation}{0}
\label{sec-5}
\subsection{Definition}
We denote by $\mathcal S_p$ the Schur class of $\mathbb C^{p\times p}$-valued Schur multipliers, meaning that the operator $M_S$ of multiplication by
$S\in\mathcal S_p$ on the left is a contraction from the Hardy space $(\mathbf H _2(\mathbb D))^p$.
In this section we consider the case of $\mathbb C^{p\times p}$-valued Schur multipliers for the space $\mathbf H_2(\mathbb K)$.
The operator-valued case, where now $S$ takes values in $\mathcal L(\mathfrak H,\mathbb C^{p\times p})$ for some Hilbert space $\mathfrak H$ will be
considered in Section \ref{sec-structure}.

\begin{theorem}
Let $S(Z)=\sum_{j=0}^\infty Z^jS_j$. Then the kernel
\begin{equation}
K_S(Z,W)=\sum_{k=0}^\infty Z^k(I_p-S(Z)S(W)^*)W^{*k}
\label{KS}
\end{equation}
is positive definite on $\mathbb K$ if and only if $s(z):=S(zI_p)$ belongs to the Schur class $\mathcal S_p$.
\label{T:Schur}
\end{theorem}
\begin{proof}
The ``only if" part: the positivity of the kernel \eqref{KS} means that for any choice of 
matrices $P_1,\ldots,P_n\in\mathbb K$, the following matrix is positive semidefinite: 
\begin{equation}
\left(\sum_{k=0}^\infty P_i^k(I_p-S(P_i)S(P_j)^*)P_j^{*k}\right)_{i,j=1,\ldots n}\geq 0.
\label{KS1}
\end{equation}
Letting $P_i=z_i I_p$ ($|z_i|<1$) gives 
$$
\left(\sum_{k=0}^\infty z_i^k(I_p-s(z_i)s(s_j)^*)\overline{z_j}^k\right)_{i,j=1,\ldots ,n}
=\left(\frac{I_p-s(z_i)s(s_j)^*}{1-z_i\overline{z_j}}\right)_{i,j=1,\ldots, n}\geq 0.
$$
Therefore, the kernel 
$$
K_s(z,w)=\frac{I_p-s(z)s(w)^*}{1-z \overline{w}}
$$
 is positive on $\mathbb D$ and hence $s$ is a $\mathbb C^{p\times p}$-valued Schur-class function.

\smallskip

Conversely, let us assume that the function $s(z)=\sum_{j=0}^\infty S_j z^j$ belongs to $\mathcal S_p$. Then it admits a coisometric
(observable) realization 
$$
s(z)=S_0+zC(I_{\mathcal X}-zA)^{-1}B=S_0+\sum_{j=1}^\infty z^j CA^{j-1}B
$$
with the state space $\mathcal X$ (the de Branges-Rovnyak space $\mathcal H(K_s)$, for example). Therefore,
$$
S_j=CA^{j-1}B\quad\mbox{for all}\quad j\ge 1.
$$
Then 
$$
S(Z)=S_0+\sum_{j=1}^\infty Z^j CA^{j-1}B=S_0+ZC\star (I_{\mathcal X}-ZA)^{-\star}B.
$$
Then the kernel \eqref{KS} is positive definite on $\mathbb K$, by Theorem \ref{T:5.2}, below.
\end{proof}

Having in view the quaternionic setting we recall the following result, which complements Theorem \ref{T:Schur}. For a proof, see \cite{donoghue}.

\begin{theorem} 
  Let $s$ be a matrix-valued function {\it defined} on a subset of the open unit disk, having an accumulation point in the open unit disk
  (as opposed to on the unit circle). Assume that the kernel $K_s(z,w)$ is positive-definite on $\Omega$. Then $s$ is the restriction to $\Omega$ of a 
  uniquely defined function analytic and contractive in the open unit disk.
  \label{T-schr-q}
\end{theorem}

Theorem \ref{T-schr-q} is used in the theory of slice-holomorphic function taking $\Omega$ to be some open subinterval of $(-1,1)$.
See \cite{zbMATH06658818}.


\begin{theorem}
    Let $S$ be a function from $\mathbb K$ into $\mathbb C^{p\times p}$. The following are equivalent:\\
    $(1)$ The function
    \begin{equation}
      \label{ksab}
K_S(A,B)=\sum_{n=0}^\infty A^n(I_p-S(A)S(B)^*)B^{*n}
\end{equation}
is positive definite on $\mathbb K$.\smallskip

$(2)$ The function $S$ is a power series: $S(Z)=\sum_{n=0}^\infty Z^nS_n$ and the operator of $\star$ multiplication by $S$ on the
left is a contraction from  $\mathbf H_2(\mathbb K)$ into itself.
\label{ksab2}
\end{theorem}

\begin{proof}
  Assume first the kernel \eqref{ksab} positive definite in $\mathbb K$. Setting $A=B$, and since $K_S(A,A)\ge 0$, we get that
  \[
\sum_{n=0}^\infty A^nS(A)S(A)^*A^{*n}\le \sum_{n=0}^\infty A^nA^{*n}<\infty,\quad A\in\mathbb K,
\]
and so the function $\sum_{n=0}^\infty Z^nS(A)^*A^{*n}C$ belongs to $\mathbf H_2(\mathbb K)$ for every $C\in\mathbb C^{p\times p}$. The positivity of the kernel then implies that the linear relation spanned by the pairs
\begin{equation}
  \label{line}
  (\sum_{n=0}^\infty Z^nA^{*n}C,\sum_{n=0}^\infty Z^nS(A)^*A^{*n}C)\subset \mathbf H_2(\mathbb K)\times\mathbf H_2(\mathbb K)
\end{equation}
extends to the graph of a contraction, say  $T$. One then computes
\[
T^*(Z^mD)=S(Z)\star (Z^mD).
\]

Conversely, assume that the operator $M_S$ of $\star$ multiplication by $S$ is a contraction from  $\mathbf H_2(\mathbb K)$ into itself. We compute
$M_S^*(K(\cdot, A)C)$ for $A\in\mathbb K$ and $C\in\mathbb C^{p\times p}$:

\[
  \begin{split}
    [M_S^*(K(\cdot, A)C), Z^mD]_{\mathbf H_2(\mathbb K)}&=[K(\cdot,A,)C, Z^mS(Z)D]_{\mathbf H_2(\mathbb K)}\\
    &=[Z^mS(Z)D,K(\cdot,A,)C]^*_{\mathbf H_2(\mathbb K)}\\
    &=\left(C^*A^mS(A)^*D\right)^*\\
    &=D^*S(A)^*A^{*m}C
  \end{split}
\]
and so, by continuity,
\[
  \begin{split}
    [M_S^*(K(\cdot, A)C), K(\cdot, B)D]_{\mathbf H_2(\mathbb K)}&=\sum_{n=0}^\infty[K(\cdot,A,)C, Z^nS(Z)B^{*n}D]_{\mathbf H_2(\mathbb K)}\\
    &=\sum_{n=0}^\infty D^*B^nS(A)^*A^{*n}C
  \end{split}
\]
and so

\begin{equation}
  M_S^*(K(\cdot, A)C)=\sum_{n=0}^\infty Z^nS(A)^*A^{*n}C 
\end{equation}
and similarly
\[
  M_S^*(K(\cdot, B)D)=\sum_{n=0}^\infty Z^nS(B)^*B^{*n}D.
\]
  Thus
\[
  [M_S^*(K(\cdot, B)D), M_S^*(K(\cdot, A)C)]=C^*\sum_{n=0}^\infty A^nS(A)S(B)^*B^{*n}D.
\]
The result follows by expressing that $M_S$ is a contraction.
\end{proof}


  \begin{remark}
    \label{rk3}
  We note that the block matrix representation of $M_S$ is
\begin{equation}
      \begin{pmatrix}S_0&0&0& \cdots&\cdots\\
        S_1&S_0    &0&\ddots &\cdots\\
        S_2& S_1&S_0     &0&\cdots\\
        \ddots&\ddots&\ddots&\ddots &\vdots
        \end{pmatrix}
    \end{equation}
  This follows from the convolution formula
  \begin{equation}
    G_n=\sum_{u=0}^nS_uF_{n-u},
    \quad u=0,1,\ldots
    \end{equation}
    for the coefficients  $G_0,G_1,\ldots\in\mathbb C^{p\times p}$ of $S\star F=\sum_{n=0}^\infty
    Z^nG_n$, where\\ $F=\sum_{n=0}^\infty Z^nF_n\in\mathbf H_2(\mathbb K)$.\smallskip

By the definition of the inner product. the contractivity is expressed as
\begin{equation}
  \label{matrix-conv} 
\sum_{n=0}^\infty \left(\sum_{u=0}^nS_u^*F_{n-u}\right)^*\left(\sum_{u=0}^nS_u^*F_{n-u}\right)\le \sum_{n=0}^\infty F_n^*F_{n}.
\end{equation}
\end{remark}

\begin{remark} (analytic extension)  The positivity of the kernel \eqref{ksab} implies that $M_S$ is a contraction from $\mathbf H_2(\mathbb K)$ into itself. In particular
  $S=M_S(I_p)\in\mathbf H_2(\mathbb K)$, and is of the form $S=\sum_{n=0}^\infty Z^nS_n$. When $p=1$ a much stronger result holds
  (see \cite{donoghue}): if (for $p=1$),
  $S$ is supposed {\bf defined} on an open set, say $\Omega$, of $\mathbb D$ (or more generally a subset of $\mathbb D$ having an acculumation point in $\mathbb D$)
  and if the corresponding kernel \eqref{ksab} is positive in $\Omega$, then 
  $S$ is the restriction to $\Omega$ of a uniquely defined function $S$ analytic and contractive in the open unit disk. We do not know if there is
  a counterpart of the result for $p>1$.
\end{remark}

\begin{remark}
We note that $K_S(A,B)$ is the unique solution of the matrix equation
\begin{equation}
  \label{main-equa}
  K_S(A,B)-AK_S(A,B)B^*=I_p-S(A)S(B)^*,\quad A,B\in\mathbb K.
\end{equation}
  This ``structural identity'' is the tool needed to extend the theory of $\mathcal H(S)$ spaces (see e.g., \cite{MR3497010,MR3617311}) from the complex scalar case to the present setting. In the case of quaternions  a similar equation and analysis hold; see \cite{acs1}. But an important
  difference is that $q\overline{q}\in[0,\infty)$ for a quaternion $q$ with conjugate $\overline{q}$, corresponding here to \eqref{eq-quart}.
\end{remark}

Note that, with $M_A^r$ being the operator of multiplication {\sl on the right} by $A$, i.e.
\[
M_A^r(Z^mB)=Z^mBA,\quad A,B\in\mathbb C^{p\times p},\,\, n\in\mathbb N,
\]
we have:
\begin{eqnarray}
  M_SM_Z&=&M_ZM_S\\
  M_SM_A^r&=&M_A^rM_S
\end{eqnarray}

We now present a version of the Bochner-Chandrasekharan theorem (see \cite[Theorem 72, p. 144]{boch_chan}) in our present setting.

\begin{theorem}
  Let $T$ be linear and contractive from $\mathbf H_2(\mathbb K)$ into itself, and satisfying
  \begin{eqnarray}
    \label{TZZT}
    TM_Z&=&M_ZT\\
    \label{TZZT2}
  TM_A^r&=&M_A^rT,\quad \forall A\in\mathbb C^{p\times p}.
  \end{eqnarray}
  Then $T$ is a Schur multiplier.
\end{theorem}

\begin{proof}
  Let $A\in\mathbb C^{p\times p}$ and $n\in\mathbb N$. We have from \eqref{TZZT}
  \begin{equation}
    TM_{Z^n}A=Z^nT(A),\quad 
  \end{equation}
  and
  \begin{equation}
    T(A)=T(I_pA)=T(I_p)A.
    \end{equation}
  Thus
  \[
  \begin{split}
    T\left(\sum_{n=0}^\infty Z^nF_n\right)&=\sum_{n=0}^\infty Z^nT(F_n) \quad\hspace{8.5mm}({\rm using \,\, \eqref{TZZT}})\\
  &=\sum_{n=0}^\infty Z^nT(I_p)F_n\quad\hspace{6mm}({\rm using \,\, \eqref{TZZT2}})\\
&=T(I_n)\star \left(\sum_{n=0}^\infty Z^nF_n\right)\end{split}
\]
\end{proof}

\begin{theorem}
Let $S=\sum_{n=0}^\infty Z^nS_n$ be a Schur multiplier, with $S_n\in\mathbb C^{p\times p}$. Then, $\widetilde{S}=\sum_{n=0}^\infty Z^nS_n^*$ is also a Schur multiplier
\end{theorem}

\begin{proof} We proceed in a number of steps.\smallskip
  
  STEP 1: {\sl Assume $S$ be a Schur multiplier. Then the function $S(zI_p)$ is analytic and contractive in the open unit disk.}\smallskip

  Setting $A=zI_p$ and $B=wI_p$ in \eqref{main-equa} with $z,w\in\mathbb C$ we get that the kernel
  \[
\sum_{n=0}^\infty z^n(I_p-S(zI_p)S(I_pw)^*)\overline{w}^n=\frac{I_p-S(zI_p)S(I_pw)^*}{1-z\overline{w}}
   \]
is positive definite in the open unit disk (since elements of the form $zI_p$ with $z$ in the open unit disk $\mathbb D$ are inside $\mathbb K$) and so
  series $S(zI_p)=\sum_{n=0}^\infty z^nS_n$ converges in the open unit disk $\mathbb D$ and is contractive there (and so is a $\mathbb C^{p\times p}$-valued Schur function of the open unit disk).\\

  STEP 2: {\sl The operator of multiplication by the function
\begin{equation}
      S(\overline{z})^*=\sum_{n=0}^\infty z^nS_n^*
    \end{equation}
    on the right  is a contraction from $\mathbf H_2(\mathbb D)^{p\times p}$ into itself (i.e. the function
    $S(\overline{z})^*$ is still a $\mathbb C^{p\times p}$-valued Schur function of the open unit disk).}\smallskip

  Using the integral representation of the inner product in $\mathbf H_2(\mathbb D)^{p\times p}$ and since
  \[
    (S(e^{-it})^*)^*S(e^{-it})^*\le I_p,\quad a.e.\,\,\mbox{\rm on the unit circle $\mathbb T$,}
    \]
  we get
  \[
\frac{1}{2\pi} \int_0^{2\pi}(f(e^{it}))^*(S(e^{-it})^*)^*S(e^{-it})^*f(e^{it})dt\le\frac{1}{2\pi} \int_0^{2\pi}f(e^{it})^*f(e^{it})dt
\]
where $f\in \mathbf H_2(\mathbb D)^{p\times p}$ and where the inequality is between matrices, and so
 \[
{\rm Tr}\, \left(\frac{1}{2\pi} \int_0^{2\pi}f(e^{it})^*(S(e^{-it})^*)^*S(e^{-it})^*f(e^{it}dt\right)\le{\rm Tr}\,\left(\frac{1}{2\pi} \int_0^{2\pi}f(e^{it})^*f(e^{it})dt\right).
\]
which expresses the asserted contractivity.\\

STEP 3: {\sl  The matrix representation of the operator multiplication by the function $S(\overline{z})^*$ on the right  has block matrix representation the 
block-Toeplitz operator with block matrix representation}
    \begin{equation}
      \begin{pmatrix}S_0^*&0&0& \cdots&\cdots\\
        S_1^*&S_0^*     &0&\ddots &\cdots\\
        S_2^*& S_1^*&S_0^*     &0&\cdots\\
        \ddots&\ddots&\ddots&\ddots &\vdots
        \end{pmatrix}
    \end{equation}


  As in Remark \ref{rk3} this follows from the convolution formula
  \begin{equation*}
    G_n=\sum_{u=0}^nS_u^*F_{n-u},\quad u=0,1,\ldots
    \end{equation*}
    for the coefficients  $G_0,G_1,\ldots\in\mathbb C^{p\times p}$ of $S(\overline{z})^* f(z)=\sum_{n=0}^\infty
    z^nG_n$, where $f(z)=\sum_{n=0}^\infty z^nF_n\in\mathbf H_2(\mathbb C^{p\times p}))$.\\

    STEP 4: {\sl $\widetilde{S}$ is Schur multiplier in $\mathbf H_2(\mathbb K)$}\smallskip

    This just follows from \eqref{matrix-conv} with the $S_n^*$ in lieu of the $S_n$ since the contractivity is expressed in the same way on the level of the coefficients in both cases.
\end{proof}



\subsection{Leech Theorem}
The main results in this section are based on discussions of the authors
with Professor Bolotnikov.
We thank Professor Bolotnikov for allowing us to include the material presented in this subsection.
\begin{theorem}
\label{T:5.2}
Let $\mathcal X$ be a Hilbert space and  let us assume that the operator
\begin{equation}
U=\begin{pmatrix}A & B \\ C & D\end{pmatrix}: \; \begin{pmatrix}\mathcal X \\ \mathbb C^p\end{pmatrix}\to
\begin{pmatrix}\mathcal X \\ \mathbb C^p\end{pmatrix}
\label{5.1}
\end{equation}
is a contraction. Then the power series
\begin{equation}
S(Z)=D+Z C\star(I_{\mathcal X}-ZA)^{-\star}B=D+\sum_{k=0}^\infty Z^{k+1}CA^k B
\label{5.2}
\end{equation}
is a Schur multiplier. Moreover, the kernel $K_S$ \eqref{KS} can be expressed as
\begin{equation}
K_S(Z,W)=\Gamma(Z)\Gamma(W)^*+\sum_{k=0}^\infty \Lambda_k(Z) (I-UU^*)\Lambda_k(W)^*
\label{5.31}
\end{equation}
where
\begin{equation}\label{gammap}
\Gamma(Z)=\sum_{k=0}^\infty Z^kCA^k\quad\mbox{and}\quad
\Lambda_k(Z)=Z^k\begin{pmatrix}Z\Gamma(Z) & I_p\end{pmatrix}.
\end{equation}
\end{theorem}
\begin{proof}
Since $U$ is a contraction, $A$ is bounded by $1$ in norm, and the $\mathcal L(\mathcal X,\mathbb C^p)$--valued power series 
${\displaystyle\sum_{k=0}^\infty Z^nCA^n}$ converges in norm for each $Z\in\mathbb K$ . It follows from \eqref{5.2} and \eqref{gammap} that
$$
\Gamma(Z)=C+Z\Gamma(Z)A, \quad S(Z)=D+ZC\Gamma(Z)B
$$
and subsequently, 
$$
\Lambda_k(Z)U=Z^k\begin{pmatrix}\Gamma(Z) & S(Z)\end{pmatrix}.
$$
Therefore,
\begin{align*}
&\Gamma(Z)\Gamma(W)^*+\sum_{k=0}^\infty \Lambda_k(Z) (I-UU^*) \Lambda_k(W)^*\\
&=\Gamma(Z)\Gamma(W)^*+\sum_{k=0}^\infty Z^k\left(I_p+Z\Gamma(Z)\Gamma(W)^*W^*\right)W^{*k}\\
&\quad-\sum_{k=0}^\infty Z^k\left(\Gamma(Z)\Gamma(W)^*+S(Z)S(W)^*\right)W^{*k}\\
&=\sum_{k=0}^\infty Z^k\left(I_p-S(Z)S(W)^*\right) W^{*k}=K_S(Z,W)
\end{align*}
which confirms \eqref{5.31}. Since the kernel on the right side of \eqref{5.31} is positive on
$\mathbb K$, the proof is complete.
\end{proof}
Our next result is Leech's factorization theorem in the present setting.
\begin{theorem}
Given two power series $P(Z)$ and $Q(Z)$, the following are equivalent:
\begin{enumerate}
\item There is a Schur multiplier $S(Z)$ such that
\begin{equation}
Q(Z)=(P\star S)(Z)\quad\mbox{for all}\quad Z\in\mathbb K.
\label{leecha}
\end{equation}
\item The kernel
\begin{equation}
\label{posab1}
K_{P,Q}(Z,W)= \sum_{k=0}^\infty
Z^k(P(Z)P(W)^*-Q(Z)Q(W)^*)W^{*k}
\end{equation}
is positive on $\mathbb K$.
\end{enumerate}
\label{T:leech}
\end{theorem}
\begin{proof} One direction is easy: if \eqref{leecha} holds for some 
$$
S(Z)=\sum_{k=0}^\infty Z^kS_k,
$$
then for every $Z, \, W\in\mathbb K$ and any $n\ge 0$, we have $Z^nP(z)=P(Z)\star Z^n$  and 
$$
Z^nQ(Z)=Z^n(P\star S)(Z)=\sum_{k=0}^\infty Z^{n+k}P(Z)S_k=P(Z)\star( Z^nS(Z)),
$$
and subsequently,
\begin{align}
K_{P,Q}(Z,W)=& \sum_{n=0}^\infty Z^n(P(Z)P(W)^*-Q(Z)Q(W)^*)W^{*n}\notag\\
=&\sum_{n=0}^\infty P(Z)\star\left(Z^nW^{*n}-Z^nS(Z)S(W)^*W^{*n}\right)\star_r P(W)^*\notag\\
=&P(Z)\star K_S(Z,W)\star_r P(W)^*\notag,
\end{align}
and the latter kernel is positive, as is seen on translating the above equality on the level of the coefficients; see \cite[Proposition 5.3 p. 855]{MR3192300} for a similar
argument for quaternionic kernels.\smallskip

Conversely, let us assume that the kernel \eqref{posab1} is positive on $\mathbb K$. Then 
(see e.g. \cite{MR647140})
it admits a Kolmogorov factorization, i.e., there exists a Hilbert space ${\mathcal X}$
and an ${\mathcal L}(\mathcal X, \mathbb C^p)$-valued power series $H(Z)$ converging weakly for all $Z\in\mathbb K$
such that 
\begin{equation}
\label{posab2}
K_{P,Q}(Z,W)= H(Z)H(W)^*\quad\mbox{for all}\quad Z,W\in\mathbb K.
\end{equation}
Combining the latter identity with \eqref{posab1} we conclude that
\begin{equation}
H(Z)H(Z)^*-ZH(Z)H(W)^*W^*=P(Z)P(W)^*-Q(Z)Q(W)^*
\label{5.4a}
\end{equation}
for all $Z,W\in\mathbb K$. The latter identity tells us that the linear map
\begin{equation}
V: \; \begin{pmatrix} H(W)^*W^*x \\
P(W)^*x\end{pmatrix}\mapsto
\begin{pmatrix}H(W)^*x\\ Q(W)^*x\end{pmatrix}
\label{defV}
\end{equation}
extends by linearity and continuity to an isometry (still denoted by $V$) from
$$
\mathcal D_V=\bigvee_{W\in\mathbb K, x\in\mathbb C^p}
\begin{pmatrix}H(W)^*W^* x \\
P(W)^* x\end{pmatrix}\subset \mathcal X\oplus\mathbb C^p
\quad \mbox{onto}\quad \mathcal R_V=\bigvee_{W\in\mathbb K, x\in\mathbb C^p}
\begin{pmatrix} H(W)^* x \\ Q(W)^*x\end{pmatrix}\subset \mathcal X\oplus\mathbb C^p.
$$
Let us extend $V$ to a contraction
$$
\widehat V=\begin{pmatrix} A^* & C^* \\ B^* & D^*\end{pmatrix}: \;
\begin{pmatrix}\mathcal X \\ \mathbb C^p\end{pmatrix}\to
\begin{pmatrix}\mathcal X \\ \mathbb C^p\end{pmatrix}.
$$
Then we have
$$
\begin{pmatrix} A^* & C^* \\ B^* & D^*\end{pmatrix}\begin{pmatrix}H(W)^*W^* x \\
P(W)^* x \end{pmatrix}= \begin{pmatrix} H(W)^* x \\ Q(W)^*x\end{pmatrix},
$$
from which it follows that 
\begin{align}
A^*H(W)^*W^*x +C^*P(W)^* x&= H(W)^*x ,\label{5.6a}\\
B^*H(W)^*W^*x+D^*P(W)^*x&= Q(W)^*x.\label{5.6b}
\end{align}
Since $A$ is a contraction, we recover $H(W)^*x$ from \eqref{5.6a} as
$$
H(W)^*x=\sum_{n=0}^\infty A^{*n}C^*P(W)^*W^*x.
$$
Substituting the latter representation into \eqref{5.6b} gives
$$
Q(W)x=D^*P(W)^*x+B^*\sum_{n=0}^\infty A^{*n}C^*P(W)^*W^{*(n+1)} x.
$$
Taking adjoints and using the arbitrariness of $x\in\mathbb C^p$ we get
\begin{align}
Q(W)&=P(W)D + \sum_{n=0}^\infty W^{n+1}P(W)CA^nB\notag\\
&=P(W)\star \left(D+\sum_{n=0}^\infty W^{n+1}CA^nB\right)\label{posab3}
\end{align}
holding for all $W\in\mathbb K$. The formula
$$
S(W)=D+\sum_{n=0}^\infty W^{n+1}CA^nB
$$
defines a Schur multiplier by Theorem \ref{T:5.2}. On the other hand, equality \eqref{posab3} means that
\eqref{leecha} is in force.
\end{proof}

\subsection{The coisometric realization}
\label{5.3}
We  associate to every Schur multiplier a co-isometric realization with state space the reproducing kernel Hilbert space module $\mathcal H(S)$ with reproducing kernel $K_S(A,B)$.
The complex-variable case was first developed, using complementation theory by de Branges and Rovnyak in \cite{dbr2}.
Here, we use the theory of linear relations as applied in \cite{adrs} for the complex-variable setting and in \cite[\S5]{AP1} for the time-varying case.\smallskip

We consider the linear relation $\mathcal R\subset (\mathcal H(S)\oplus\mathbb C^{p\times p})\times ({\mathcal H(S))}\oplus\mathbb C^{p\times p})$ spanned by the pairs
\[
\left(  \begin{pmatrix} K_S(\cdot, A)A^*G\\ H\end{pmatrix}\, ,\,\begin{pmatrix} (K_S(\cdot , A)-K_S(\cdot, 0))G+K_S(\cdot, 0)H\\ (S(A)^*-S(0)^*)A^*G+S(0)^*H\end{pmatrix}\right)
\]
when $A$ runs through $\mathbb K$ and $G,H$ run through $\mathbb C^{p\times p}$.

\begin{proposition}
$\mathcal R$ is densely defined and isometric, and thus extends to the graph of an everywhere defined isometry from ${\mathbf H_2(\mathbb K)}\oplus\mathbb C^{p\times p}$ into itself.
\end{proposition}

\begin{proof}
  Let
\[
\left(  \begin{pmatrix} K_S(\cdot, B)B^*E\\ F\end{pmatrix}\, ,\,\begin{pmatrix} (K_S(\cdot , B)-K_S(\cdot, 0))E+K_S(\cdot, 0)H\\ (S(B)^*-S(0)^*)E+S(0)^*F\end{pmatrix}\right)
\]
be another element of $\mathcal R$, with $B\in\mathbb K$ and $E,F\in\mathbb C^{p\times p}$. We want to show that
\[
\begin{split}
  \langle K_S(\cdot, B)B^*E, K_S(\cdot, A)A^*G\rangle_S+{\rm Tr}\,H^*F&=\\
  &\hspace{-6cm}=\langle (K_S(\cdot , B)-K_S(\cdot, 0))E+K_S(\cdot, 0)F,(K_S(\cdot , A)-K_S(\cdot, 0))G+K_S(\cdot, 0)H\rangle_S+\\
  &\hspace{-5.5cm}+{\rm Tr}((S(B)^*-S(0)^*)E+S(0)^*F)^*((S(A)^*-S(0)^*)G+S(0)^*H).
  \end{split}
\]
Considering the $\mathbb C^{p\times p}$-valued forms associated to the inner products, the above  equality can be rewritten in the form
\[
  G^*\square_1E+G^*\square_2F+H^*\square _3E+H^*\square _4F=0,
  \]
  for appropriate expressions $\square_j$, $j=1,2,3,4$ which we now show to be equal to $0$. We have
  \[
    \begin{split}
      \square_1&=-AK_S(A,B)B^*+K(A,B)-K_S(A,0)-K_S(0,B)+K_S(0,0)\\
      &\hspace{5mm}+(S(A)-S(0))(S(B)^*-S(0)^*)\\
&=I_p-S(A)S(B)^*-(I_p-S(A)S(0)^*)-(I_p-S(0)S(B)^*)+(I_p-S(0)S(0)^*)+\\
      &\hspace{5mm}+S(A)S(B)^*-S(A)S(0)^*-S(0)S(B)^*+S(0)S(0)^*\\
      &=0,
      \end{split}
    \]
    where we have used \eqref{main-equa} to go from the first equality to the second one.\smallskip

    Similarly
    \[
      \begin{split}
        \square_2&=K_S(A,0)-K_S(0,0)+(S(A)-S(0))S(0)^*\\
        &=I_p-S(A)S(0)^*-I_p+S(0)S(0)^*+S(A)S(0)^*-S(0)S(0)^*\\
        &=0,
        \end{split}
      \]
and
          \[
      \begin{split}
        \square_3&=K_S(0,B)-K_S(0,0)+S(0)(S(B)^*-S(0)^*)\\
        &=I_p-S(0)S(B)^*-I_p+S(0)S(0)^*+S(0)S(B)^*-S(0)S(0)^*\\
        &=0.
        \end{split}
      \]
As for $\square_4$ one has:
                \[
        \square_4=-I_p+K_S(0,0)+S(0)S(0)^*=0.
      \]
      To conclude we prove that the domain of $\mathcal R$ is dense in  $\mathcal H(S)\oplus\mathbb C^{p\times p}$.
      Any element $\begin{pmatrix} F\\ C\end{pmatrix}\in\mathcal H(S)\oplus \mathbb C^{p\times p}$ orthogonal to the domain of $\mathcal R$ will satisfy
      \[
        {\rm Tr} (G^*AF(A)+H^*C)=0.
      \]
      Letting $G=0$ and $H$ run through $\mathbb C^{p\times p}$ we get $C=0$. We thus have        $ {\rm Tr} (G^*AF(A)=0$ for all $G\in\mathbb C^{p\times p}$,
      and so $AF(A)=0$ and so $aF(aI_p)=0$ for all $a$ near the origin and so that $F=0$.
\end{proof}

We write the isometry whose graph extends $\mathcal R$ in the form
\[
  \begin{pmatrix}\mathscr T&\mathscr F\\
    \mathscr G&\mathscr H\end{pmatrix}^*,
\]
that is
\begin{eqnarray}
  \label{6-1}
  \mathscr T^*\left( K_S(\cdot,A)A^*G\right)&=&(K_S(\cdot , A)-K_S(\cdot, 0))G\\
  \label{6-2}
  \mathscr F^* \left(K_S(\cdot,A)A^*G\right)&=&(S( A)^*-S(0)^*)G\\
  \label{6-3}
  \mathscr G^*H&=&K_S(\cdot, 0)H\\
  \mathscr H^*H&=&S(0)^*H.
  \end{eqnarray}

  \begin{proposition}
Let $F\in\mathcal H(S)$, $A\in\mathbb K$ and $G\in\mathbb C^{p\times p}$.  It holds that
\begin{eqnarray}
  \label{6-5}
A(\mathscr T F)(A)&=&F(A)-F(0)\\
A( \mathscr F G)(A)&=&(S(A)-S(0))G.\\        
\mathscr G F&=&F(0),\\
    \mathscr H H&=&S(0)H.
                  \end{eqnarray}
                \end{proposition}

                \begin{proof}

Apply both sides of \eqref{6-1} to $F$ to obtain
  \begin{equation}
    \label{7-1}
    \langle F ,\mathscr T^*\left( K_S(\cdot,A)A^*G\right)\rangle_S=\langle F,K_S(\cdot , A)-K_S(\cdot, 0))G\rangle_S,
  \end{equation}
  that is
  \[
    G^*A(\mathscr T F)(A)=G^*(F(A)-F(0)).
    \]
                  \end{proof}

                  \begin{theorem}
                    Let $S$ be a Schur multiplier. In the above notation we have
                    \begin{equation}
                      \label{tyui}
                      S(A)G=\mathscr H G+\sum_{n=0}^\infty A^{n+1}\mathscr G\mathscr T^n\mathscr F(G),\quad G\in\mathbb C^{p\times p}.
\end{equation}
                    \end{theorem}

                    \begin{proof}
                      Let $C_A$ denote evaluation at $A\in\mathbb K$. Equation \eqref{6-5} gives
                      \[
                        C_A=C_0+AC_A\mathscr T
                      \]
                      from which we get
                      \[
                        C_A=\sum_{n=0}^\infty A^nC_0\mathscr T^{n}
                      \]
                      which converges in $\mathcal H(S)$ in the operator norm.
                      Applying to $F=\mathscr F G$ we obtain
                      \[
(\mathscr F G(A))=\sum_{n=0}^\infty A^nC_0\mathscr T^{n+1}\mathscr F G,
\]
and so
\[
  A(\mathscr F G(A))=\sum_{n=0}^\infty A^{n+1}C_0\mathscr T^{n+1}\mathscr F G,
\]
i.e.
\[
  S(A)G-S(0)G=\sum_{n=0}^\infty A^{n+1}C_0\mathscr T^{n+1}\mathscr F G.
  \]
This concludes the proof since $C_0=\mathscr G$ and  $S(0)=\mathscr H$.
\end{proof}

\subsection{Operator ranges and complementation}
\label{5.4}
In this section we give an another proof of the realization theorem for Schur multipliers using complementation theory and the theory of operator ranges. We refer to  \cite[\S4]{fw} and and \cite{sarason94} for the latter. Similar arguments for the time-varying setting can be found in \cite{AP1}, and in
\cite{a-et-al} for the quaternionic setting. We first recall a well-known result on operator ranges. We provide a proof for completeness.

\begin{theorem}
  \label{thran}
  Let $\mathfrak H$ be a Hilbert space and let $A$ be a bounded positive operator from $\mathfrak H$ into itself. Let $\pi$ denote the orthogonal projection
  onto $\ker A$.
  Then the space ${\rm ran}\, \sqrt{A}$ endowed with the norm (called the range norm)
  \begin{equation}
    \label{norm}
  \|\sqrt{A}h\|_{{\rm ran}\,\sqrt{A}}=\|(I_{\mathfrak H}-\pi)h\|_{\mathfrak H},\quad h\in\mathfrak H,
  \end{equation}
  is a Hilbert space. Furthermore,
  \begin{equation}
    \|Ah\|_{{\rm ran}\,\sqrt{A}}=\|\sqrt{A}h\|_{\mathfrak H},\quad h\in\mathfrak H,
  \end{equation}
  and the range of $A$ is dense in the range of $\sqrt{A}$ in this norm.
  Finally  in the associated inner product it holds that
  \begin{equation}
    \label{polar}
\langle \sqrt{A}h, \sqrt{A}g\rangle_{{\rm ran}\,\sqrt{A}}=\langle (I_{\mathfrak H}-\pi)h,g\rangle_{\mathfrak H}
  \end{equation}
  \begin{equation}
\langle Ah, Ag\rangle_{{\rm ran}\,\sqrt{A}}=\langle Ah,g\rangle_{\mathfrak H}
  \end{equation}
  and
  \begin{equation}
    \label{inner2}
    \langle Ah, \sqrt{A}g\rangle_{{\rm ran}\,\sqrt{A}}=\langle \sqrt{A}h,g\rangle_{\mathfrak H}
  \end{equation}
  for $ h,g\in\mathfrak H$.
    \end{theorem}

\begin{proof}
  That \eqref{norm} is indeed a norm follows from
  \[
  \sqrt{A}h=0\quad\iff \pi h=h.
  \]
  Let now a Cauchy sequence $(\sqrt{A}h_n)$ in ${\rm ran}\,\sqrt{A}$. Then $(I_{\mathfrak H}-\pi)h_n$ is a Cauchy sequence in $\mathfrak H$ converging to an
  element, say $k$. We have by continuity of $\pi$ that $k=(I_{\mathfrak H}-\pi)k$ and so $\sqrt{A}h_n)$ has limit $\sqrt{A}k$ and ${\rm ran}\,
  \sqrt{A}$ is closed in the   range norm. Formula \eqref{polar} follows from \eqref{norm} by polarization. The last two formulas follow from
  \[
  A(I_{\mathfrak H}-\pi)=A.
  \]
\end{proof}

We now recall the operator range characterization from \cite{fw} which we will use.
Note that \cite{fw} consider one Hilbert space and we consider two possibly different Hilbert spaces, but the proof is the same for the latter case.
Theorem \ref{thth} below is  \cite[Theorem 4.1 p. 275]{fw} with $\sqrt{I_{\mathfrak H}-TT^*}$ rather than $T$. We present a proof since the passage
from one to the other involves  {\it a priori} polar representations. We give a direct proof, explicitly adapted from \cite{fw}. See also
\cite{a-et-al} where it is presented in the quaternionic setting.
\begin{theorem}
  \label{thth}
  Let $T$ be a contraction from the Hilbert space $\mathfrak G$ into the Hilbert space $\mathfrak H$. Then, $f\in{\rm ran}\,\sqrt{I_{\mathfrak H}- TT^*}$
  if and only if
  \begin{equation}
    \label{marseille}
    \sup_{g\in\mathfrak G}\left(\|f+Tg\|_{\mathfrak G}^2-\|g\|^2_{\mathfrak G}\right)<\infty.
\end{equation}
  \end{theorem}

\begin{proof}
  As mentioned above, the proof is directly adapted from the proof of \cite[Theorem 4.1 p. 275]{fw} with $\sqrt{I_{\mathfrak H}-TT^*}$ rather than $T$.
  We have for $g\in\mathfrak G$ and $h\in\mathfrak H$
  \[
    \begin{split}
      \|\sqrt{I_{\mathfrak H}-TT^*}h+Tg\|_{\mathfrak H}^2-\|g\|^2_{\mathfrak G}&=\\
      &\hspace{-2cm}
      =\|\sqrt{I_{\mathfrak H}-TT^*}h\|^2_{\mathfrak H}
+\|Tg\|_{\mathfrak H}^2-\|g\|^2+2{\rm Re}\, \langle \sqrt{I_{\mathfrak H}-TT^*}h,Tg\rangle_{\mathfrak H}      
-\|g\|^2_{\mathfrak G}\\
&\hspace{-2cm}=\|h\|_{\mathfrak H}^2-\|T^*h\|_{\mathfrak G}^2+2{\rm Re}\, \langle T^*h,\sqrt{I_{\mathfrak G}-T^*T}g\rangle_{\mathfrak G}      -\|\sqrt{I_{\mathfrak G}-T^*T}
g\|^2_{\mathfrak G}\\
&\hspace{-2cm}=\|h\|_{\mathfrak H}^2-\|T^*h-\sqrt{I_{\mathfrak G}-T^*T}
g\|^2_{\mathfrak G}\\
&\hspace{-2cm}\le \|h\|^2_{\mathfrak H}
      \end{split}
    \]

    Conversely,   let $f\in\mathfrak H$ such that the supremum in \eqref{marseille} is finite, and denote it by $K$. The choice $g=0$ in \eqref{marseille}
    gives $\|u\|^2_{\mathfrak H}\le K$. We can write
\[
\|f+Tg\|^2_{\mathfrak H}\le K+\|g\|^2_{\mathfrak G},\quad\forall g\in\mathfrak G
\]
and so
\begin{equation}
  \label{rty123}
2{\rm Re}\, \langle f,Tg\rangle_{\mathfrak H}\le K-\|u\|^2_{\mathfrak H}+\|g\|^2_{\mathfrak G}-\|Tg\|^2_{\mathfrak H}
\end{equation}
If $\langle f,Tg\rangle_{\mathfrak H}=0$, this inequality is trivial.  If $\langle f,Tg\rangle_{\mathfrak H}\not=0$ and replace $g$ by $tge^{i\theta}$
where $t\in\mathbb R$ and $\theta\in
\mathbb R$ is chosen such that
\[
e^{-i\theta}\langle f,Tg\rangle_{\mathfrak H}=|\langle f,Tg\rangle_{\mathfrak H}|.
\]
Using the previous equality for every $g$ we rewrite \eqref{rty123} as
\[
2t|\langle f,Tg_1\rangle_{\mathfrak H}|\le K-\|f\|_{\mathfrak H}^2+t^2\|\sqrt{I_{\mathfrak G}-T^*T}g_1\|^2_{\mathfrak G},\quad \forall t\in\mathbb R,
\,\,\forall g_1\in\mathfrak H.
\]
It follows that
\[
  0  \le t^2\|\sqrt{I_{\mathfrak G}-T^*T}g_1\|^2_{\mathfrak G}-2t|\langle f,Tg_1\rangle_{\mathfrak H}|+K-\|f\|_{\mathfrak H}^2,
  \quad \forall t\in\mathbb R,\,\,\forall g_1\in\mathfrak H.
\]
Thus
\[
  |\langle f,Tg_1\rangle_{\mathfrak H}|\le \sqrt{K-\|f\|^2_{\mathcal H}}
  \left(\|\sqrt{I_{\mathfrak G}-T^*T}g_1\|_{\mathfrak G}\right).  \]
  By Theorem \ref{thran} we have
  \[
\|\sqrt{I_{\mathfrak G}-T^*T}g_1\|_{\mathfrak G}=\|(I_{\mathfrak G}-T^*T)g_1\|_{{\rm ran}\, (I_{\mathfrak G}-T^*T)}
\]
and since the range of $I_{\mathfrak G}-T^*T$ is dense in the range of $\sqrt{I_{\mathfrak G}-T^*T}$ in the range norm, 
the map $g_1\mapsto \langle Tg_1, u\rangle_{\mathfrak H}$ is continuous on the range of $\sqrt{I_{\mathfrak G}-T^*T}$.
By Riesz representation theorem  for continuous linear functionals on a Hilbert space  there exists $r\in\mathfrak H$
  such that
  \[
  \begin{split}
    \langle Tg_1, u\rangle_{\mathfrak H}&=\langle (I_{\mathfrak G}-T^*T)    g_1,\sqrt{I_{\mathfrak G}-T^*T} r\rangle_{{\rm ran}\, \sqrt{I_{\mathfrak G}-T^*T}}\\
&=    \langle \sqrt{I_{\mathfrak G}-T^*T}g_1, r\rangle_{\mathfrak G}
\end{split}
     \]
as follows from \eqref{inner2}.  It follows that $T^*f=\sqrt{I_{\mathfrak G}-T^*T}r$. So
  \[
  \begin{split}
    f&=f-TT^*f+TT^*f\\
    &=(I_{\mathfrak H}-TT^*)f+T\sqrt{I_{\mathfrak G}-T^*T}r\\
    &=(I_{\mathfrak H}-TT^*)f+\sqrt{I_{\mathfrak G}-TT^*}r
    \end{split}
    \]
    belongs to ${\rm ran}\, \sqrt{I_{\mathfrak H}-TT^*}$.
\end{proof}

We now apply the previous results to our present setting.

\begin{proposition}
  It holds that
  \[
\mathcal H(S)=\left\{ G\in {\mathbf H_2(\mathbb K)}\,\, ;\sup_{H\in{\mathbf H_2(\mathbb K)}}\|G+M_SH\|^2_{\mathbf H_2(\mathbb K)}-\|H\|^2_{\mathbf H_2(\mathbb K)}<\infty\right\}
\]
and the above supremum is then the norm of $G$ in $\mathcal H(S)$.
\end{proposition}

\begin{proof}
This is Theorem \ref{thth} with $\mathfrak H=\mathbf H_2(\mathbb K)$ and $T=M_S$.
  \end{proof}

\begin{proposition}
  It holds that
  \begin{equation}
    \mathcal H(S)={\rm ran}\,\sqrt{I_{\mathbf H_2(\mathbb K)}      -M_SM_S^*}
\end{equation}
with the operator range norm
\[
\|\sqrt{I_{\mathbf H_2(\mathbb K)}-M_SM_S^*}F\|=\|(I_{\mathbf H_2(\mathbb K)}-\pi)F\|_{\mathbf H_2(\mathbb K)}
\]
where $\pi$ is the orthogonal projection onto $\ker \sqrt{I_{\mathbf H_2(\mathbb K)}-M_SM_S^*}$.
\end{proposition}

\begin{proof}
  This is Theorem \ref{thran} with $\mathfrak H=\mathbf H_2(\mathbb K)$ and $T=M_S$.
\end{proof}

\begin{theorem}
  \begin{equation}
    \mathbf H_2(\mathbb K)={\rm ran}\, M_S+{\rm ran}\,\sqrt{I_{\mathbf H_2(\mathbb K)}-M_SM_S^*}
  \end{equation}
in the sense of complementation.
  \end{theorem}
\subsection{Another approach to the co-isometric realization}
We study the co-isometric realization of a Schur multiplier using the results of the previous section, that is, using complementation theory. See \cite{AP1}
for similar computations in the time-varying setting. Recall that $R_0$ was defined in \eqref{r00}.

\begin{proposition}
  Let $F\in\mathcal H(S)$.Then, $R_0F\in\mathcal H(S)$ and
  \begin{equation}
    \|R_0F\|_{\mathcal H(S)}^2\le\|F\|_{\mathcal H(S)}^2-\|F(0)^*F(0)\| 
  \end{equation}
  Let $C\in\mathbb C^{p\times p}$. Then, $R_0SC\in\mathcal H(S)$ and
    \begin{equation}
      \|R_0(SC)\|_{\mathcal H(S)}^2\le {\rm Tr}\, (C(I_p-S(0)S(0)^*)C^*).
      \end{equation}
  \end{proposition}

  \begin{proof}
    Let $G\in{\mathbf H_2(\mathbb K)}$.
\[
  \begin{split}  
    \|R_0F+M_SG\|_{\mathbf H_2(\mathbb K)}^2-\|G\|_{\mathbf H_2(\mathbb K)}^2&\\
    &\hspace{-3cm}=\|M_Z(R_0F+M_SG)\|_{\mathbf H_2(\mathbb K)}^2-\|G\|_{\mathbf H_2(\mathbb K)}^2\\
    &\hspace{-3cm}=\|F-F_0+M_SZG)\|_{\mathbf H_2(\mathbb K)}^2-\|G\|_{\mathbf H_2(\mathbb K)}^2\\
    &\hspace{-3cm}=\|-F_0\|_{\mathbf H_2(\mathbb K)}^2-2{\rm re}\langle F_0,F+M_SZG\rangle_{\mathbf H_2(\mathbb K)}+\\
    &\hspace{-2.5cm}+   \underbrace{ \|F+M_SZG\|_{\mathbf H_2(\mathbb K)}^2-\|ZG\|_{\mathbf H_2(\mathbb K)}^2}_{\le \|F\|_{\mathcal H(S)}^2}\\
    &\hspace{-3cm}\le\|F\|^2_{\mathcal H(S)}-\|F_0F_0^*\|_{\mathbb C^{p\times p}.}
\end{split}
\]
 Similarly   
  \[
  \begin{split}
    \|R_0(SC)+M_SG\|_{\mathbf H_2(\mathbb K)}^2-\|G\|_{\mathbf H_2(\mathbb K)}^2&\\
    &\hspace{-3cm}=\|M_Z(R_0(SC+M_SG)\|_{\mathbf H_2(\mathbb K)}^2-\|G\|_{\mathbf H_2(\mathbb K)}^2\\
    &\hspace{-3cm}=\|SC-S_0C+M_S(ZG)\|^2_{\mathbf H_2(\mathbb K)}-\|G\|_{\mathbf H_2(\mathbb K)}^2\\
    &\hspace{-3cm}=\|-S_0C+M_S(C+ZG)\|^2_{\mathbf H_2(\mathbb K)}-\|G\|_{\mathbf H_2(\mathbb K)}^2\\
    &\hspace{-3cm}=\|-S_0C\|^2_{\mathbf H_2(\mathbb K)}+2{\rm Re}\,\langle -S_0C, M_S(C+ZG)\rangle_{\mathbf H_2(\mathbb K)}^2+\\
    &\hspace{-2.5cm}+\|M_S(C+ZG)\|^2_{\mathbf H_2(\mathbb K)}-\|M_ZG\|_{\mathbf H_2(\mathbb K)}^2\\
    &\hspace{-3cm}=\|S_0C\|^2_{\mathbf H_2(\mathbb K)}-2\|S_0C\|^2 _{\mathbf H_2(\mathbb K)}+\|M_S(C+ZG)\|^2_{\mathbf H_2(\mathbb K)}-\|G\|_{\mathbf H_2(\mathbb K)}^2\\
    &\hspace{-3cm}
    \le-\|S_0C\|^2_{\mathbf H_2(\mathbb K)}+\|C+ZG\|^2_{\mathbf H_2(\mathbb K)}-\|G\|_{\mathbf H_2(\mathbb K)}^2\\
    &\hspace{-3cm}=\|C\|^2_{\mathbf H_2(\mathbb K)}-\|S_0C\|^2_{\mathbf H_2(\mathbb K)}.
  \end{split}
  \]
  \end{proof}
\subsection{A structure theorem}
\label{sec-structure}
An important aspect of the theory of de Branges and de Branges-Rovnyak spaces (see \cite{dbjfa,dbr1,dbr2}) is the weakening of isometric inclusion to
contractive inclusion, and the associated notion of complementation, which replaces orthogonal sum. The results involve matrix-valued, or
more generally operator-valued functions. For example (see \cite[p. 24]{ad3}) the function
\[
s(z)=\begin{pmatrix}c_1s_1(z)&c_2s_2(z)&\cdots&c_Ns_N(z)\end{pmatrix}
  \]
  where $s_1,\ldots, s_N$ are inner functions (for instance finite Blaschke products) and $c_1,c_2,\ldots, c_N$ are complex numbers such that
  \[
    \sum_{n=1}^N|c_n|^2=1
  \]
  is such that
  \[
\frac{1-s(z)s(w)^*}{1-z\overline{w}}=\sum_{n=1}^N |c_n|^2\frac{1-s_n(z)\overline{s_n(w)}}{1-z\overline{w}},
    \]
    and the associated reproducing kernel Hilbert space will, in general, be only contractively included in the Hardy space $\mathbf H_2(\mathbb D)$.
    \smallskip
    
We also remark that in the theory of reproducing kernel Hilbert spaces, complementation and contractive inclusion
correspond to the older results on the reproducing kernel Hilbert space associated to a sum of positive definite functions; see \cite[p. 353]{aron},
\cite{MR0023243}. In Section \ref{Cara} this problem is avoided by assuming that one starts with a power series to begin with.\smallskip

In this section we replace the inequality
\[
  \|R_0F\|_{\mathfrak M}^2=\|F\|^2_{\mathfrak M}-\|F(0)\|^2
\]
which characterizes isometric inclusion in $\mathbf H_2(\mathbb K)$, by the inequality
\begin{equation}
  \label{yuiop}
  \|R_0F\|_{\mathfrak M}^2\le\|F\|^2_{\mathfrak M}-\|F(0)\|^2.
\end{equation}
So, we wish to study the structure of $R_0$-invariant subspaces contractively included in $\mathbf H_2(\mathbb K)$. We follow the arguments in the proof
of \cite[Theorem 3.1.2 p. 85]{adrs}, suitably adapted to the present situation. First note that \eqref{yuiop} can be rewritten as
\[
  R_0^*R_0+C^*C\le I_{\mathfrak M}.
\]

So
\[
  \begin{pmatrix}R_0\\ C\end{pmatrix}^* \begin{pmatrix}R_0\\ C\end{pmatrix}\le I_{\mathfrak M}.
      \]
Since the adjoint of a Hilbert space contraction is a Hilbert space contraction we can write
\[
  \begin{pmatrix}R_0\\ C\end{pmatrix}  \begin{pmatrix}R_0\\ C\end{pmatrix}^*\le I_{{\mathfrak M}\oplus \mathbb C^p}.
\]
Let $\mathfrak H=\mathfrak M\oplus \mathbb C^{p\times p}$  and let $B\in\mathcal L(\mathfrak H,\mathfrak M)$ and $D\in
\mathcal L(\mathfrak H,\mathbb C^{p\times p})$ be defined by
\[
\begin{pmatrix}B\\ D\end{pmatrix}=\sqrt{I_{{\mathfrak M}\oplus\mathbb C^{p\times p}}-  \begin{pmatrix}R_0\\ C\end{pmatrix}  \begin{pmatrix}R_0\\ C\end{pmatrix}^*}
\]
Then,
\[
M=  \begin{pmatrix}
    R_0&G\\
    C&D\end{pmatrix}\,\,\,\mathfrak M\oplus \mathfrak H\,\,\longrightarrow\,\,\,\mathfrak H\oplus \mathbb C^{p\times p}
\]
is co-isometric. We define
\[
  S=D+\sum_{n=1}^\infty Z^nCR_0^{n-1}G.
\]
Note that $S$ is $\mathcal L(\mathfrak H,\mathbb C^{p\times p})$-valued and that $S(A)$ makes sense for all $A\in\mathbb K$ as a converging series
in the operator norm topology

\begin{remark}
When $Z=\lambda I_p$ the function $S$ coincides with the characteristic function of the colligation $M$; see e.g. \cite[p. 16]{adrs}.
\end{remark}

Following \cite[p. 38]{MR2240272}, we associate to $S$ the (possibly unbounded and not everywhere defined)
multiplication operator $M_S$ now from $\ell_2(\mathbb N_0,\mathfrak H)$ into $\mathfrak H_2(\mathbb K)$ by
\begin{equation}
M_S(h)=\sum_{n=0}^\infty Z^n\left(Dh_n+\sum_{j=0}^{n-1}CR_0^{n-1-j}Gh_j\right),\quad h=(h_j)_{j=0}^\infty\in \ell_2(\mathbb N_0,\mathfrak H).
  \end{equation}

\begin{theorem}
  $M_S$ is a contraction from $\ell_2(\mathbb N_0,\mathfrak H)$ into $\mathbf H_2(\mathbb K)$ if and only if the kernel \eqref{ksab} (now computed for the current
  operator-valued $S$) is positive definite in $\mathbb K$.
\end{theorem}

\begin{proof}
We have (using the co-isometry of  $M$)
\[
  \begin{split}
    I_p-S(A)S(B)^*&=\underbrace{I_p-DD^*}_{CC^*}-\sum_{n=1}^\infty A^nCR_0^{n-1}GD^*-\\
    &\hspace{5mm}-\sum_{m=1}^\infty \underbrace{DG^*}_{-CR_0^*}R_0^{*(m-1)}C^*B^{*m}-\sum_{n,m=1}^\infty A^nCR_0^{n-1}GG^*R_0^{*(m-1)}C^*B^{*m}\\
    &=CC^*+\sum_{n=1}^\infty A^nCR_0^{n-1}R_0C^*+\sum_{m=1}^\infty CR_0^*R_0^{*(m-1)}C^*B^{*m}-\\
    &\hspace{5mm}    -\sum_{n,m=1}^\infty A^nCR_0^{n-1}(I-R_0R_0^*)R_0^{*(m-1)}C^*B^{*m}\\
    &=\underbrace{CC^*}_{\stackrel{{\rm def.}}{=}\mathbf 1} +\underbrace{\sum_{n=1}^\infty A^nCR_0^nC^*}_{\stackrel{{\rm def.}}{=}\mathbf 2}
    +\underbrace{\sum_{m=1}^\infty CR_0^{*m}C^*B^{*m}}_{\stackrel{{\rm def.}}{=}\mathbf 3}+\\
    &\hspace{5mm}  +\underbrace{\sum_{n,m=1}^\infty A^nCR_0^{n}R_0^{*m}C^*B^{*m}}_{\stackrel{{\rm def.}}{=}\mathbf 4}  -
    \underbrace{\sum_{n,m=1}^\infty A^nCR_0^{n-1}R_0^{*(m-1)}C^*B^{*m}}_{\stackrel{{\rm def.}}{=}\mathbf 5}\\
    &=\mathbf 1+\mathbf 2+\mathbf 3+\mathbf 4-\mathbf 5
  \end{split}
  \]
so that
\begin{equation}
  \label{above}
I_p-S(A)S(B)^*=\underbrace{\sum_{n,m=0}^\infty A^nCR_0^nR_0^{*m}C^*B^{*m}}_{\mathbf 1+\mathbf 2+\mathbf 3+\mathbf 4}
- \underbrace{\sum_{n.m=0}^\infty A^{n+1}CR_0^{n}R_0^{*m}C^*B^{*(m+1)}}_{\mathbf 5}
\end{equation}

  Define now
  \[
K(A,B)=\left(\sum_{n=0}^\infty A^nCR_0^n\right)\left(\sum_{m=0}^\infty B^mCR_0^m\right)^*=\sum_{n,m=0}^\infty A^nCR_0^nR_0^{*m}C^*B^{*m}
\]
Then $K(A,B)$ is positive definite and
\[
  \begin{split}
    K(A,B)-AK(A,B)B^*&=\sum_{n,m=0}^\infty A^nCR_0^nR_0^{*m}C^*B^{*m}-\sum_{n,m=0}^\infty A^{n+1}CR_0^{n}R_0^{*m}C^*B^{*(m+1)}\\
    &=I-S(A)S(B)^*
    \end{split}
  \]
  by \eqref{above}, and so
  \begin{equation}
    \label{kasab3}
K(A,B)=\sum_{n=0}^\infty A^n(I_p-S(A)S(B)^*)B^{*n},
\end{equation}
which is \eqref{ksab}, but now with $S$ operator-valued. To prove that $M_S$ is a contraction we adapt the proof of Theorem \ref{ksab2} as follows. Recall that $S(A)^*$ is a bounded operator from $\mathbb C^{p\times p}$ into $\mathfrak H$.
The linear relation \eqref{line} is now a linear subspace of $\mathbf H_2(\mathbb K)\times \ell_2(\mathbb N_0,\mathfrak H)$ spanned by the pairs
\[
(\sum_{n=0}^\infty Z^nA^{*n}C, (S(A)^*A^{*n}C)_{n=0}^\infty).
\]
The positivity of \eqref{kasab3} implies that this linear relation extends to the graph of a contraction, whose adjoint is $M_S$.\\

The converse statement is a direct computation.
\end{proof}

\begin{remark}
For $A=aI_p$ and $B=bI_p$ the function $K(A,B)$ reduces to formulas given in \cite[Theorem 2.1.2 p. 44 and p. 97]{adrs}.
  \end{remark}

  \begin{remark}
    Besides the scalar case $A=aI_p$ with $a\in\mathbb C$,  we do not know if and when the positivity of the kernel \eqref{kasab3} implies that
    $S(A)$ is a contraction. The counterpart of this question in the quaternionic setting has
    a negative answer in the matrix-valued case; see \cite[(62.38) p. 1767]{zbMATH06526253}. We adapt the example from the latter publication to the present setting to find a
    counterexample for multipliers of     $\mathbf H_2(\mathbb K)$. We consider $(\mathbf H_2(\mathbb K))^2$. Multipliers are elements $S=(S_{uv})_{u,v=1,2}$ of
    $(\mathbf H_2(\mathbb K))^{2\times 2}$ such that
  \[
  \sum_{n=0}^\infty A^n\left(I_{2p}-\begin{pmatrix}S_{11}(A)&S_{12}(A)\\ S_{21}(A)& S_{22}(A)\end{pmatrix}
  \begin{pmatrix}(S_{11}(B))^*&(S_{21}(B))^*\\ (S_{12}(B))^*&
      (S_{22}(B))^*\end{pmatrix}\right)B^{*n}
  \]
  is positive definite in $\mathbb K$. Let  $p\ge 2$ and  $J\in\mathbb C^{p\times p}$ be such that $J^*=-J$ and $J^2=-I_p$ (for instance, if $p=2$,
  $J=i\begin{pmatrix}1&0\\ 0&-1  \end{pmatrix}$). Let
  \[
S=\frac{1}{\sqrt{2}}\begin{pmatrix}Z&J\\ZJ&I\end{pmatrix}
\]
(compare with \cite[(62.38) p. 1767]{zbMATH06526253}), and let
\[
  F=\sum_{n=0}^\infty Z^n\begin{pmatrix}A_n\\ B_n\end{pmatrix}\in(\mathbf H_2(\mathbb K))^2.
  \]
Then,
\[
  \begin{split}
    (\sqrt{2}S)\star F&=\sum_{n=0}^\infty\begin{pmatrix}Z^{n+1}A_n+Z^nJB_n\\ Z^{n+1}JA_n+Z^nB_n\end{pmatrix}\\
    &=\begin{pmatrix}JB_0\\ B_0\end{pmatrix}+\sum_{k=0}^\infty Z^{2k+1}\begin{pmatrix}A_{2k}+JB_{2k+1}\\ JA_{2k}+B_{2k+1}\end{pmatrix}+
    \sum_{k=1}^\infty Z^{2k}\begin{pmatrix}A_{2k-1}+JB_{2k}\\ JA_{2k-1}+B_{2k}\end{pmatrix}.
    \end{split}
  \]
  Since
  \[
    \begin{split}
      (A_{2k}+JB_{2k+1})^*(A_{2k}+JB_{2k+1})+(JA_{2k}+B_{2k+1})^*(JA_{2k}+B_{2k+1})&=\\
      &\hspace{-5cm}=A_{2k}^*A_{2k}+B_{2k+1}^*B_{2k+1}-B_{2k+1}^*JA_{2k}+B_{2k+1}^*JA_{2k}\\
      &\hspace{-5cm}=A_{2k}^*A_{2k}+B_{2k+1}^*B_{2k+1}
  \end{split}
\]
and similarly
  \[
    \begin{split}
      (A_{2k-1}+JB_{2k})^*(A_{2k-1}+JB_{2k})+(JA_{2k-1}+B_{2k})^*(JA_{2k-1}+B_{2k})&=\\
      &\hspace{-5cm}=A_{2k-1}^*A_{2k-1}+B_{2k}^*B_{2k}
  \end{split}
\]
we have
  \[
  [S\star F ,S\star F]_{(\mathbf H_2(\mathbb K))^2}=[F, F]^2_{(\mathbf H_2(\mathbb K))^2},
\]
and so $S$ is a Schur multiplier. But 
\[
  I_{2p}-S(A)S(A)^*=\frac{1}{2}\begin{pmatrix}I_p-AA^*&AJA^*-J\\ AJA^*-J&I_p-AA^*\end{pmatrix}
\]
which is not positive for $A$ unitary such that $AJA^*\not =J$ ( we need $p\ge 2$ for ensure this; take for instance
$A=\frac{1}{\sqrt{2}}\begin{pmatrix}1&1\\1&-1\end{pmatrix}$ and $J$ as above.
    \end{remark}

    \section{Carath\'eodory multipliers}
\label{Cara}
    \setcounter{equation}{0}
    Closely related to Schur functions are Carath\'eodory functions, that is, functions analytic in the open unit disk and with a positive real part there. In
    \cite{herglotz}, see also the collection of papers \cite{hspnw}, Herglotz gave an integral representation fo such functions. A function $\varphi$ is a Carath\'eodory function if and only if it can be written as
    \begin{equation}
      \label{rep-hergl}
      \varphi(z)=im+\int_0^{2\pi}\frac{e^{is}+z}{e^{is}-z}d\mu(s),\quad z\in\mathbb D,
    \end{equation}
    where $m\in\mathbb R$ and where $\mu$ is a positive measure on $[0,2\pi)$. Note that \eqref{rep-hergl} can be rewritten as
    \begin{equation}
      \label{rep-hergl-1}
      \varphi(z)=im+\int_0^{2\pi}d\mu(s)+2\sum_{n=1}^\infty z^n t_n
    \end{equation}
    where
    \begin{equation}
t_n=\int_{0}^{2\pi}e^{-ins}d\mu(s),\quad n=0,1,\ldots
      \end{equation}
      is the moment sequence associated to $d\mu$.\smallskip
      
    In terms of kernels, a function $\varphi$ {\sl defined} on a subset of the
    open unit disk which possesses an accumulation point in $\mathbb D$ is the restriction to $\Omega$ of a (uniquely defined) Carath\'eodory
    function if and only if the kernel
    \begin{equation}
      \label{v-kernel}
      \frac{\varphi(z)+\overline{\varphi(w)}}{1-z\overline{w}}
      \end{equation}
      is positive definite in $\Omega$. Note that \eqref{v-kernel} can be rewritten as
      \[
      \frac{\varphi(z)+\overline{\varphi(w)}}{1-z\overline{w}}=\sum_{n=0}^\infty z^n (\varphi(z)+\overline{\varphi(w)})\overline{w}^n.
    \]
    \begin{remark}
We note that $\varphi$ need not be bounded in modulus in the open unit disk. Then (and only then), the operator of multiplication by $\varphi$ will not be a bounded operator from the Hardy space into itself.
\end{remark}

        The previous discussion motivates the following definition:

        \begin{definition}
          A $\mathbb C^{p\times p}$-valued function $\Phi$ defined in $\mathbb K$ is called a Carath\'eodory multiplier if the kernel
          \begin{equation}
     K_\Phi(A,B)=     \sum_{n=0}^\infty A^n(\Phi(A)+\Phi(B)^*)B^{*n}
            \label{Phi-123456}
\end{equation}
is positive definite on $\mathbb K$.
\end{definition}
        In this section we prove a realization theorem similar to \eqref{rep-hergl} for such functions, in two different ways:

\begin{enumerate}
  \item The first approach reduces the study to the complex setting, and does not require that the function $\Phi(A)$, or more generally,  that the operator $M_\Phi$ of $\star$-multiplication on the left be
    a bounded operator from $\mathbf H_2(\mathbb K)$ into itself.

    \item 
    The second approach will require this latter hypothesis. Then, $\Phi$ is a power series in $Z$ with matrix coefficients on the right, and converging in $\mathbb K$, i.e. $\Phi=\sum_{n=0}^\infty Z^n\Phi_n$.
    As just mentionned above, In the classical setting, this hypothesis is not necessary.
    \end{enumerate}

Rather than stating the theorem and proving it afterwards, we here prefer to go the other way around, and begin with a discussion and results which lead to the result. The result itself is presented in Theorem
\ref{5678!!!!} below. So let us start from a function $\Phi$ for which the kernel \eqref{Phi-123456} is positive definite in $\mathbb K$, and set
$\Psi(a)=\Phi(aI_p)$ with $a\in\mathbb D$. The kernel
\begin{equation}
\frac{\Psi(a)+\Psi(b)^*}{1-a\overline{b}}
  \end{equation}
  is positive definite in the open unit disk. By the matrix version of Herglotz representation theorem (see e.g. \cite[Theorem 4.5 p. 23]{MR48:904}
  for the operator-valued version), we can write
  \[
\Psi(a)=iX+\int_{0}^{2\pi}\frac{e^{it}+z}{e^{it}-z}dM(t)
    \]
    where $X\in\mathbb C^{p\times p}$ is self-adjoint and where $M$ is now a $\mathbb C^{p\times p}$-valued positive measure on $[0,2\pi)$. This formula can be rewritten as
    \begin{equation}
\Psi(a)=iX+T_0+2\sum_{n=1}^\infty a^nT_n
\end{equation}
where $(T_n)_{n\in\mathbb N_0}$ is the moment sequence associated to $dM$. We claim that (recall that $\Phi$ is assumed to be a power series in $Z$ with matrix coefficients on the right, converging in $\mathbb K$)
\begin{equation}
  \label{phi-789}
  \Phi(A)=iX+T_0+2\sum_{n=1}^\infty A^nT_n,\quad A\in\mathbb K.
\end{equation}
Indeed, let $\widetilde{\Phi}$ denote the right hand side of \eqref{phi-789}. Both $\widetilde{\Phi}$ and $\Phi$ coincide on the matrices $aI_p$, $a\in\mathbb D$, and this restriction completely determines the coefficients $T_n$. Thus we have proved one direction in the following result. The converse direction
is easily proved and will be omitted.
\begin{theorem}
  Let $\Phi=\sum_{n=0}^\infty Z^n\Phi_n$ be a power series in $Z$ with matrix coefficients on the right, converging in $\mathbb K$ and assume $\Phi_0=\Phi_0^*$. Then the kernel
  \[
    \sum_{n=0}^\infty A^n(\Phi(A)+\Phi(B)^*)B^{*n}
  \]
  is positive definite in $\mathbb K$ if and only if the sequence $(\Phi_n)_{n\in\mathbb N_0}$ is the moment sequence of a positive $\mathbb C^{p\times p}$-valued measure on $[0,2\pi)$.
\label{5678!!!!}
\end{theorem}

We note that the method used here, different from the one we used to characterize Schur multiplier, is not intrinsic, in the sense that we do not use the reproducing kernel Hilbert space with reproducing kernel \eqref{Phi-123456}. Denoting this space by $\mathcal L(\Phi)$, one can also characterize $\Phi$ in
terms of the associated backward-shift realization, assuming the operator $M_\Phi$ bounded from $\mathbf H_2(\mathbb K)$ into itself. 

\begin{theorem}
  \label{567890}
  Assume that $M_\Phi$ is bounded from $\mathbf H_2(\mathbb K)$ into itself, and let $\mathcal L(\Phi)$ denote the reproducing kernel space with reproducing kernel $K_\Phi$. Then
  ${\rm Re}\, M_\Phi\ge 0$ and
  \begin{equation}
    \mathcal L(\Phi)={\rm ran}\,\sqrt{M_\Phi+M_\Phi^*}
  \end{equation}
  with the operator range norm.
\end{theorem}

The proof follows the arguments of Section \ref{5.4} and will be omitted.
The key in the proof is the formula
\begin{equation}
M_\Phi^*((I-ZB)^{-\star}E)=\sum_{n=0}^\infty Z^n\Phi(B)^*B^{*n}E,\quad B\in\mathbb K,\quad E\in\mathbb C^{p\times p},
  \end{equation}
valid since $M_\phi$ is bounded, and so for $A,B\in\mathbb K$ and $E,H\in\mathbb C^{p\times p}$,
\begin{equation}
  \begin{split}
    [ (M_\Phi+M_\Phi^*)((I-ZB)^{-\star}E),(I-ZA)^{-\star}H]_{\mathcal L(\Phi)}=\sum_{n=0}^\infty A^n(\Phi(A)+(\Phi(B))^*)B^{*n}.
  \end{split}
  \end{equation}

\smallskip

The second realization theorem is now presented and proved.

\begin{theorem}
  Assume that the operator $M_\phi$ is bounded from $\mathbf H_2(\mathbb K)$ into itself.  Then $R_0$ is a co-isometry from $\mathcal L(\Phi)$ into itself.
  Furthermore, $\Phi$ can be written as $\Phi=\sum_{n=0}^\infty Z^n\Phi_n$
\begin{eqnarray}
\label{phi-coeff-1}
  {\rm Re}\, \Phi_0&=&2C_0C_0^*\\
  \label{phi-coeff}
  \Phi_n&=&C_0R_0^nC_0^*,\quad n=1,\ldots
\end{eqnarray}
where $C_0$ is the evaluation at $0$.
\end{theorem}
    \begin{proof}
      We follow the approach from \cite[pp. 708-709]{atv1}, suitably adapted to the present setting. We first note that $K_\Phi$ in \eqref{Phi-123456}  satisfies the equation
\begin{equation}
  \label{structure-equa}
K_\Phi(A,B)-AK_\Phi(A,B)B^*=\Phi(A)+(\Phi(B))^*,\quad A,B\in\mathbb K.
\end{equation}
      We then define in $\mathcal L(\Phi)\times \mathcal L(\Phi)$ the linear relation  $\mathscr R$ spanned by the pairs (which we write as a row rather than column, in opposition to the notation in Section \ref{5.3})
      \begin{equation}
        (K_\Phi(\cdot, B)BE,K_\Phi(\cdot, B)E-K_\Phi(\cdot,0)E),\quad B\in\mathbb K,\quad E\in\mathbb C^{p\times p},
      \end{equation}
      and divide the rest of the proof into steps.\\

      STEP 1: {\sl $\mathscr R$ is isometric.}\smallskip

    Indeed, let $A\in\mathbb K$ and $H\in\mathbb C^{p\times p}$.
    We have on the one hand
    \begin{equation}
      \label{poiuy}
[ K_\Phi(\cdot,B)BE,K_\phi(\cdot, A)AH]_{\mathcal L(\Phi)}=H^*AK_\Phi(A,B)B^*E,
\end{equation}
and on the other hand,
\[
\begin{split}
  [K_\Phi(\cdot, B)E-K_\Phi(\cdot,0)E),  K_\Phi(\cdot, A)H-K_\Phi(\cdot,0)H]_{\mathcal L(\Phi)}&=\\
  &\hspace{-7cm}=H^*K_\Phi(A,B)E-H^*K_\Phi(A,0)E-H^*K_\Phi(0,B)E+H^*K_\Phi(0,0)E\\
  &\hspace{-7cm}=H^*K_\Phi(A,B)E-H^*(\Phi(A)+\Phi(0))E-H^*(\Phi(0)+(\Phi(B))^*)E+\\
  &\hspace{-6.6cm}+H^*(\Phi(0)+(\Phi(0))^*)E\\
  &\hspace{-7cm}=H^*K_\Phi(A,B)E-H^*(\Phi(A)E-H^*(\Phi(B))^*E\\
  &\hspace{-7cm}=H^*\left\{K_\Phi(A,B)E-\Phi(A)-(\Phi(B))^*\right\}E.
\end{split}
  \]
  By \eqref{structure-equa},
  \[
    H^*\left\{K_\Phi(A,B)E-\Phi(A)-(\Phi(B))^*\right\}E=H^*A^*K_\Phi(A,B)BE,
  \]
  and hence the isometry property holds.\\

  STEP 2: {\sl $\mathscr R$ has a dense domain.}\smallskip

  Indeed, let $F=\sum_{n=0}^\infty Z^nF_n$ be orthogonal to the domain of $\mathscr R$. Then for every $B$ and $E$ as above, $E^*BF(B)=0$. Taking $B=bI_p$ with $b\in\mathbb D$, this leads to $F_n=0$ for $n=0,1,\ldots$ and so $F=0$.\\

  It follows from the first two steps that $\mathscr R$ extends to the graph of an everywhere defined isometry, say $T$, which we compute in STEP 3.\\
  
  STEP 3: {\sl It holds that $T^*=R_0$, and so $\mathcal L(\Phi)$ is $R_0$-invariant.}\smallskip
  
  Let $F\in\mathcal L(\Phi)$, and $B,E$ as above. We have
  \[
    \begin{split}
      E^*BF(B)&=[TF,K_\Phi(\cdot,B)B^*E]_{\mathcal L(\Phi)}\\
      &=[F,K_\Phi(\cdot, B)E-K_\Phi(\cdot,0)E]_{\mathcal L(\Phi)}\\
      &=E^*(F(B)-F(0)),
      \end{split}
    \]
    so that
    \[
      B(TF(B))=F(B)-F(0).
    \]
    Since $\mathcal L(\Phi)\subset\mathbf H_2(\mathbb K)$ (see Theorem \ref{567890}) we know that $TF$ is a power series with coefficents on the right it follows that $TF=R_0$.\\

    STEP 4: {\sl We prove \eqref{phi-coeff}.}\smallskip

    We first note that $\Phi$ can be written as $\Phi=\sum_{n=0}^\infty Z^n\Phi_n$ since $K_\Phi(\cdot, 0)E\in\mathcal L(\Phi)$ for every
    $E\in\mathbb C^{p\times p}$. We note that
    \[
K_\Phi(\cdot, 0)=\Phi+(\Phi(0))^* \in\mathcal L(\Phi)
    \]
    and so
  \[
    C_0C_0^*=2{\rm Re}\, \Phi_0,
  \]
which is \eqref{phi-coeff-1}. Furthermore,     
\[
  \begin{split}
    C_0(K(\cdot,0)E)&=(\Phi_0+\Phi_0^*)E\\
    C_0(R_0(K(\cdot ,0)E)&=C_0(\Phi_1+Z\Phi_2+\cdots)=\Phi_1E\\
    &\hspace{1.7mm}\vdots\\
    C_0(R_0^n(K(\cdot ,0)E)&=C_0((\Phi_n+Z\Phi_n+\cdots)E)=\Phi_nE\\
        &\hspace{1.7mm}\vdots
  \end{split}
  \]
\end{proof}

\begin{remark}
  \eqref{phi-coeff} implies that, with $Z=zI_p$,
  \begin{equation}
    \begin{split}
      \Phi E&=      \sum_{n=0}^\infty z^n\Phi_nE\\
      &=    ( i {\rm Im}\,\Phi(0))E+({\rm Re}\, \Phi_0)E+\sum_{n=1}^\infty z^nC_0R_0^{n-1}C_0^*E\\
      &=   (  i {\rm Im}\,\Phi(0))E+\frac{C_0C_0^*E}{2}+\sum_{n=1}^\infty z^nC_0R_0^{n-1}C_0^*E\\
      &=   (  i {\rm Im}\,\Phi(0))E+\frac{1}{2}C_0(I_{\mathcal L(\Phi)}-zR_0)^{-1}(I_{\mathcal L(\Phi)}+zR_0)C_0^*E.
    \end{split}
  \end{equation}
  In the above expression, $z$ is a number and $zR_0$ means the multiplication by this number
  of the operator $R_0$ acting in $\mathcal L(\Phi)$.
  \end{remark}

  For papers related to this section we mention \cite{MR2003f:47021,MR1704663,MR1839830}.
  \section{Concluding remarks}
\setcounter{equation}{0}
We conclude with some remarks on possible future work.\smallskip

\subsection{Interpolation}
Problem \ref{illustration-pb} is a special case of a much more general bitangential interpolation problem (see e.g. \cite{bgr}), which can also be set in the framework of Schur multipliers. We will
consider these problems in a future publication.

\subsection{Matrix polynomials and other applications of the map \eqref{aug23bbbb}}
The map \eqref{aug23bbbb} which to $F(Z)=\sum_{n=0}^\infty Z^nF_n$ associates the matrix-valued function of a complex variable
$F(zI_p)=\sum_{n=0}^\infty z^nF_n$ allows further applications than the ones presented here. For instance,  define $F(Z)$ to be a matrix-polynomial
(of the matrix variable $Z$) if only a finite number of powers of $Z$ arise in the power series expansion of $F$. Then, $F(Z)$ is a matrix polynomial
if and only if $F(zI_p)$ is a  classical matrix-polynomials, and factorizations of $F$ in factors of matrix-polynomials coincide with
factorization of $F(zI_p)$ into classical matrix-polynomials. These factorizations, and much more, have been considered in
\cite{zbMATH03983386,zbMATH00042718}. We plan to pursue this line of study in a future work.

\subsection{Symmetries}
We begin with a motivating example.
For a general matrix $A=(a_{jk})_{j,k=1}^h\in\mathbb C^{2h\times 2h}$, define
  \begin{equation}
    A_\varphi=J_1\overline{A}J_1^*,\quad{\rm where}\quad J=\begin{pmatrix}0&I_h\\
      -I_h&0\end{pmatrix}\quad {\rm and}\quad \overline{A}=(\overline{a_{jk}})_{j,k=1}^h.
  \end{equation}
  Clearly
\[
(  \lambda I_{2h})_\varphi=\overline{\lambda}I_{2h},\quad\lambda\in\mathbb C
\]

The proof of the following lemma is easy and will be omitted.

\begin{lemma}
  \label{lemma7-1}
  We have
    \begin{equation}
      \label{asdfg}
      (A_\varphi)^*=(A^*)_\varphi\quad{\rm and} \quad     \overline{A_\varphi}=(\overline{A})_\varphi\quad{\rm and}
    \end{equation}
    
  \begin{equation}
    \label{abphi}
(AB)_{\varphi}=A_\varphi B_{\varphi}
    \end{equation}
    and 
    \[
      A\ge 0\quad\Longrightarrow\quad A_\varphi\ge0
    \]
    \end{lemma}

    Furthermore:

    \begin{lemma}
      \label{lemma7-2}
    We have
    \begin{equation}
      \label{wertyu}
      \sqrt{A_\varphi}=(\sqrt{A})_\varphi
    \end{equation}
  \end{lemma}

  \begin{proof}
    We have $A=\sqrt{A}\sqrt{A}$ and so from \eqref{abphi}
    \[
      A_\varphi=(\sqrt{A})_\varphi(\sqrt{A})_\varphi
    \]
    If $A=A_\varphi$ and since $(\sqrt{A})_\varphi\ge 0$ the uniqueness of the squareroot implies \eqref{wertyu}.
  \end{proof}

  \begin{lemma}
   If $A$ is invertible, we have
    \[
      A_{\varphi}^{-1}=(A_\varphi)^{-1}
    \]
  \end{lemma}
  \begin{proof}
    This follows from  $\overline{A}^{-1}=\overline{A^{-1}}$.
    \end{proof}
    
    Finally, a matrix $A$ satisfies $A=A_\varphi$ if and only if it is of the form

\begin{equation}
  A=\begin{pmatrix}A_1&A_2\\
  -  \overline{A_2}&\overline{A_1}
  \end{pmatrix}
\end{equation}
where $A_1$ and $A_2$ belong to $\mathbb C^{h\times h}$.\smallskip

\begin{proposition}
  Restricting $A$ and $Z$ in  the formula \eqref{UA} for the Blaschke factor $U_A(Z)$, the latter satisfies
  \[
    (U_A(Z))_\varphi=U_A(Z).
    \]
  \end{proposition}

  \begin{proof}
This follows from \eqref{wertyu}, \eqref{abphi} and \eqref{asdfg}.
    \end{proof}

    When $h=1$, $AA^*$ is a scalar matrix and we get back the classical Blaschke factor from quaternionic analysis.\\

    Similarly, let now

  \begin{equation}
    A_\varphi=J_2\overline{A}J_2^*,\quad{\rm where now }\quad J_2=\begin{pmatrix}0&I_h\\
      I_h&0\end{pmatrix}.
  \end{equation}

Then, Lemmas \ref{lemma7-1} and \ref{lemma7-2} still hold. Furthermore, a matrix $A$ now satisfies $A=A_\varphi$ if and only if it is of the form

\begin{equation}
  A=\begin{pmatrix}A_1&A_2\\
  \overline{A_2}&\overline{A_1}
  \end{pmatrix}
\end{equation}
where $A_1$ and $A_2$ belong to $\mathbb C^{h\times h}$. When $h=1$ we get back the split quaternions, and $U_A$ will be the corresponding Blaschke factor. Even when
$h=1$ we get a new formula since $AA^*$ will not be a scalar matrix.\\

\begin{remark}
  Taking $A_1$ and $A_2$ with real values we get the real realization of elements in $\mathbb C^{h\times h}$ in the first case and  $h\times h$ matrices with components
  hyperbolic numbers in the second case.
  \end{remark}

  \begin{definition}
    We call a map satisfying
      \begin{eqnarray}
                (AB)_{\varphi}&=&A_\varphi B_{\varphi}\\
                (A+B)_{\varphi}&=&A_\varphi +B_{\varphi}\\
        A\ge 0&\Longrightarrow& A_\varphi\ge0\\
(A_\varphi)^*&=&(A^*)_\varphi\\
(  \lambda I_{2h})_\varphi&=&\overline{\lambda}I_{2h},\quad\lambda\in\mathbb C
      \end{eqnarray}
      an admissible symmetry.
    \end{definition}

  For such a symmetry it follows that
    \[
      A_\varphi=(\sqrt{A})_\varphi(\sqrt{A})_\varphi
    \]
    and, when $A$ is invertible,
    \[
      A_{\varphi}^{-1}=(A_\varphi)^{-1}
    \]

    \begin{definition}
The ring of matrices for which $A=A_\varphi$ is called the associated ring and denoted by $\mathbb C_\varphi$.
    \end{definition}

    The various results presented here extend when we replace $\mathbb C^{p\times p}$ by $\mathbb C_\varphi$. This setting includes quaternions, split-quaternions and
    corresponding matrix versions.
For a related work, see \cite{MR4658650}.
  \section*{Acknowledgments} It is a pleasure to thank Professor Alain Yger for very useful comments and remarks on various versions of the
  present paper. It is a pleasure to thank Professor Bolotnikov for numerous discussions on this paper and for allowing us to include the material presented in subsection \ref{T:5.2}.
             Daniel Alpay thanks the Foster G. and Mary McGaw Professorship in Mathematical Sciences, which supported this research, and in particular the stays of Professor Cho at Chapman University.
\bibliographystyle{plain}
\def\cprime{$'$} \def\cprime{$'$} \def\cprime{$'$}
  \def\lfhook#1{\setbox0=\hbox{#1}{\ooalign{\hidewidth
  \lower1.5ex\hbox{'}\hidewidth\crcr\unhbox0}}} \def\cprime{$'$}
  \def\cprime{$'$} \def\cprime{$'$} \def\cprime{$'$} \def\cprime{$'$}
  \def\cprime{$'$}

\end{document}